%% file: arXiv.tex
\documentclass[]{article}
\usepackage{amsmath, amsthm, amssymb,bm,datetime}
\usepackage{graphicx,color,subfigure}
\usepackage{epstopdf}
\usepackage{booktabs}
\usepackage{cases}
\usepackage{cite}
\usepackage[margin=0.75in]{geometry}
\usepackage{authblk}

\include{symbols}

\begin{document}
%
\title{Accelerated Magnetic Resonance Thermometry in the Presence of Uncertainties
}
%
%
%

\author[1]{Reza~Madankan*}
\author[1]{Wolfgang~Stefan}
\author[1]{Samuel Fahrenholtz}
\author[1]{Christopher MacLellan}
\author[1]{John~Hazle}
\author[1]{Jason Stafford}
\author[2]{Jeffrey S.~Weinberg, Ganesh~Rao}
\author[1]{David~Fuentes}

\affil[1]{Department of Imaging Physics, University of Texas MD Anderson Cancer Center, Houston,
TX, 77015 USA e-mail: rmadankan@mdanderson.org}
\affil[2]{Department of Neurosurgery, University of Texas MD Anderson Cancer Center}

\date{}
\maketitle{}

\begin{abstract}
An accelerated model-based information theoretic approach is presented to
perform the task of Magnetic Resonance (MR) thermal image reconstruction from 
a limited number of observed samples on $k-$space. The key
idea of the proposed approach is to utilize information theoretic techniques to \textit{optimally detect} 
samples of $k-$space that are information rich with respect to a
model of the thermal data acquisition. These highly informative $k-$space samples are then used to  
refine the mathematical model and reconstruct the image. 
The information theoretic reconstruction is demonstrated retrospectively
in data acquired  during MR--guided Laser Induced
Thermal Therapy (MRgLITT) procedures.
{\color{black}
The approach demonstrates that locations of
high-information content with respect to a model based reconstruction  of
MR thermometry may be quantitatively identified.
The predicted locations of high-information content are sorted  and
retrospectively extracted from the fully sampled $k-$space measurements data set.
The effect of 
interactively increasing the predicted number of data points used in the
subsampled reconstruction is quantified using the 
L2-norm of the distance between the subsampled and fully sampled
reconstruction.
Performance of the proposed approach is also compared with clinically available subsampling techniques (rectilinear subsampling and variable-density Poisson disk undersampling). It is shown that the proposed subsampling scheme results in accurate reconstructions 
using small fraction of $k-$space points and
suggest that the reconstruction technique may be useful in 
improving the efficiency of the thermometry data temporal resolution.
}
\end{abstract}

\subsection*{keywords}
MRI, thermal ablation, $k-$space, inverse problems, entropy, sensor management, quantification and estimation, Pennes bioheat model.

%

\section{Introduction}
Laser Induced Thermal Therapy (LITT) is a non-ionizing and
minimally invasive treatment alternative for patients that have reached
maximum radiation dose or have surgically unresectable 
tumors~\cite{stafford2010laser,fahrenholtz2013generalised,carpentier2012mr,schwarzmaier2006mr,jethwa2012magnetic,norred2014magnetic}.
The approach is being used clinically
to rapidly treat focal diseased tissue in 
liver, prostate, bone, brain, and other types of diseases like epilepsy  \cite{rhim1999radiofrequency,chen2000prostate,anzai1995preliminary,simon2006percutaneous,curry2012mr}. 
MR thermometry is often integrated with the LITT
procedure to provide temperature estimates of the target tissue and
surrounding critical structures during the laser induced thermal therapy. 
Currently, MR thermometry is acquired primarily through
quantitative phase-sensitive techniques based on the Proton Resonance
Frequency (PRF) shift.  
The information provided by the integrated system, often called MR--guided Laser Thermal Therapy (MRgLITT), improves safety
and allows the clinician to reduce power or stop delivery if
safety thresholds are exceeded~\cite{mcnichols2004mr}. 

The idea of combining MR thermometry signals with a mathematical model to better represent the temperature field has been exploited in various works \cite{roujol2012robust,de2013extended,fuentes2012kalman}. The essence of all these works is to provide a better estimate of the temperature image based on model--data fusion in a Kalman filtration \cite{kalman1960new} framework. However, they don't provide any guidelines about the optimal $k-$space sampling.

This manuscript presents a novel model-based information theoretic approach
to increase both the accuracy and efficiency of the thermal image
acquisition and reconstruction.
The key question is that whether there exist \textit{optimal}
locations on $k-$space such that they provide \textit{maximum information}
regarding the temperature image? And if there exist, how one can develop a general
framework to find these locations on $k-$space? 

Numerous works have studied the effect of data undersampling
on temperature image reconstruction, and image reconstruction in
general
\cite{hennig1999k,willinek2002randomly,lustig2007sparse,majumdar2012compressed,todd2009temporally,chartrand2009fast,todd2010model,gaur2014accelerated,todd2012reconstruction,odeen2014sampling,song2009optimal,vogt2007high,lim20083d,candes2008introduction,vasanawala2011practical}.
Hennig \cite{hennig1999k} performed a comprehensive study on the effect of
rectilinear and non-rectilinear methods for $k-$space sampling. Ode{\'e}n \etal
\cite{odeen2014sampling} also performed a comprehensive study on different
$k-$space subsampling strategies for MR temperature monitoring in
particular. Willinek \etal \cite{willinek2002randomly} suggested a randomly
segmented central $k-$space ordering for high-spatial-resolution-contrast
enhanced MR Angiography. Lustig \etal \cite{lustig2007sparse} proposed a
sparse sampling technique based on compressed sensing theory
\cite{candes2008introduction} to reconstruct MR images from undersampled
$k-$space.
Vasanawala \cite{vasanawala2011practical} \etal developed a practical parallel imaging compressed sensing MRI by help of Poisson-disk random undersampling of the phase encode. 
 Majumdar \etal \cite{majumdar2012compressed} also applied the
compressed sensing theory for real time dynamic MR image reconstruction.
Performance of compressed sensing in non-convex optimization context for
sparse sampling of $k-$space was also studied by Chartrand
\cite{chartrand2009fast}. Stochastic sampling of $k-$space has also been
used in contrast-enhanced MR renography and time-resolved angiography
\cite{song2009optimal,vogt2007high,lim20083d}. Todd \etal
\cite{todd2009temporally,todd2012reconstruction} presented a temporally
constrained reconstruction (TCR) algorithm for reconstruction of temperature
image from undersampled $k-$space data in an optimization framework. 

Model-based reconstruction of temperature data is also recently exploited in few works. Todd \etal 
\cite{todd2010model} provided a model predictive filtering approach for temperature 
image reconstruction during the MR--guided High Intensity Focus Ultrasound (MRgHIFU). The key idea of this work is to
utilize a mathematical model to predict the spatiotemporal variations of
temperature, where its parameters are identified based on observed $k-$space
data. The identified mathematical model is then combined with a set of
under-sampled $k-$space data (by simply replacing the model predictions with
observed data) to provide an estimate of temperature image. 
Geometric techniques were used for under-sampling of $k-$space in this work.
Gaur \etal
\cite{gaur2014accelerated} proposed an accelerated MR thermometry approach by 
model-based estimation of the temperature from undersampled
$k-$space. In the heart of this approach, there exists a constrained image
model whose parameters are learned from acquired $k-$space data points.
The fitted image model is then used to reconstruct artifact-free temperature
maps. However, no scheme is provided for selection of $k-$space
samples used in finding model parameters. 

Information theoretic approaches have also been used in optimal experiment design in MRI and medical imaging in general  \cite{marseille1996bayesian,brihuega2003optimization,poot2010optimal,cercignani2006optimal}. The essence of these works lies in finding the optimal experiment parameters based on minimizing Cramer-Rao lower bound. Minimizing Cramer-Rao lower bound provides a reliable technique for finding the optimal parameters of the experiment. However, the major drawback of this technique is that the Fisher information used in Cramer-Rao lower bound can be a function of uncertain parameters. Hence, one needs to use an estimate of the parameters of interest (e.g. prior mean estimate of the parameter) to approximate the Fisher information in Cramer-Rao lower bound, resulting in possible sub-optimalities. 

While the majority of work has focused on randomly selected
or geometrically undersampled $k-$space points for image reconstruction,
herein we tackle the problem of $k-$space under-sampling and temperature image
reconstruction from a different perspective. In detail, we first apply an uncertain
mathematical model to describe the spatiotemporal evolution
of the temperature field and corresponding manifestation of $k-$space measurements. 
We then apply a general information theoretic framework to identify
the \textit{most useful} locations on $k-$space for data acquisition.
Finally, the obtained $k-$space data points are used to \textit{refine} the
mathematical model in reconstructing the temperature image over the tissue. 

The key contribution of this article is to provide a general framework to \textit{judiciously} select the
samples on $k-$space such that they result in the best estimate of the
corresponding temperature image. This has been achieved by maximizing the entropy of
model outputs or maximizing the mutual information \cite{cover2012elements}
between model uncertainties and data observations. 
\textcolor{black}{\textit{Optimally selected} observations
on $k-$space are then utilized to \textit{refine} the mathematical model.
The \textit{refined mathematical model} is then used to generate the
temperature field corresponding to the observations.}


The structure of this paper is as follows: First, mathematical description
of the problem is provided in Section \ref{sec:prob_statement}. Then, we
describe the overall picture of the proposed approach in Section
\ref{sec:methodology}. 
In section \ref{sec_uq}, we discuss the uncertainty
quantification for Pennes bioheat equation in presence of parametric uncertainty. 
The basic idea for optimal $k-$space sampling is then elucidated
in Section \ref{sec:dyn_data}. 
A brief discussion about the minimum variance
framework used for model-data fusion and model based image reconstruction is provided in Sections \ref{sec:gpc_minvariance} and \ref{sec:reconstruction}, respectively. 
In Section \ref{sec:num_sim}, performance of the
proposed algorithm is demonstrated using \textit{in silico} and \textit{in
vivo} examples. Results and discussion conclude the manuscript.

\section{Problem Statement}\label{sec:prob_statement}
Consider the problem of temperature field estimation during the MRgLITT
process. According to the Pennes bioheat equation \cite{pennes1998analysis},
spatio-temporal behavior of temperature $u(\x,t)$ at a given time $t$ and spatial location $\x$ is given by:

\begin{align}
\raggedright
\rho c \frac{\partial u}{\partial t}-\nabla\cdot({\Lambda(u,\x)}\nabla u)&+\nonumber\\
{\omega(u,\x)}c_{blood}(u-u_a)=Q_{laser}(\x,t), \quad \forall\ \x \in \Omega &\nonumber\\
Q_{laser}(\x,t)=P(t)\mu^2\frac{e^{-\mu||\x-\x_0||}}{4\pi ||\x-\x_0||}\nonumber\\
u(\x,0)=u_0(\x),\quad \x \in \Omega & \nonumber\\
u(\x,t)=u_D(\x),\quad \x \in \partial \Omega_D & \label{pennes}\\
-\Lambda(u,\x)\nabla\cdot u(\x,t)=g_N(\x),\quad \x \in \partial \Omega_N & \nonumber\\
-\Lambda(u,\x)\nabla\cdot u(\x,t)=h(u-u_{\infty}),\quad \x \in \partial \Omega_c & \nonumber
\end{align}
where, 
parameters $\rho$ and $c$ are density and specific heat of
the tissue, and $\Lambda(u,\x),\ c_{blood}$, and $\omega(u,\x)$ represent
thermal conductivity, specific heat capacity  of blood, and perfusion rate,
respectively. Note that $\x$ is an $n-$dimensional spatial domain in
general, \textit{i.e.} $\x\in \re^n$. For instance, if $n=2$, then $\Omega$ is a two
dimensional domain and $\x=(x,y)$.  

$Q_{laser}(\x,t)$ denotes the optical-thermal response to the laser source and
is modeled as the classical spherically symmetric isotropic solution to
the transport equation of light within a laser-irradiated tissue
\cite{welch1984thermal}. $P(t)$ is the laser power as a function of time. Parameter
$\x_0$ denotes the position of laser photon source, and $\mu$ is the optical
attenuation coefficient. Note that $\mu$ can be a function of tissue properties in general, i.e. $\mu=\mu(\x)$.

\textcolor{black}{The initial temperature of the tissue is denoted by $u_0(\x)$, which is
typically assumed to be $37\ ^{\circ}$C for human body. Also, $u_D(\x)$ and $g_N(\x)$ represent Dirichlet and Neumann's boundary conditions, respectively. Possible effects of blood flow on temperature are also modeled as convection $h(u-u_{\infty})$, where $h$ denotes the convection coefficient and $u_{\infty}$ represents the blood temperature. Note that $u_{\infty}=u_0$.
}

One of the key challenges while using \eq{pennes} is that the associated
parameters like optical attenuation coefficient $\mu(\x)$ are often unknown,
and their values are different from patient to patient. Hence, it is a
reasonable assumption to treat these parameters as \textit{uncertain
parameters} with some \textit{a priori} information. Here we will 
assume the optical attenuation coefficient to be uncertain with some prior
probability density function $p(\mu)$. Further, tissue heterogenieties are
incorporated using finite spatial variations of the attenuation coefficient:

\begin{equation}\label{conduction}
\mu(\x)=\sum\limits_{i=1}^d \mu_i\ U(\x-\Omega_i)
\end{equation}
where, $d$ is the total number of tissue types and $\mu_i$ is the
corresponding attenuation coefficient for the $i^{th}$ tissue. Note that
each of the $\mu_i$ is assumed to be uncertain. $\Omega_i$
represents the corresponding sub-domain of the $i^{th}$ tissue on spatial
domain $\Omega$. Note that $\Omega_i$'s are mutually disjoint, i.e.
\[\bigcup_{i=1}^d \Omega_i = \Omega,\quad \Omega_i\cap \Omega_j=\varnothing
\]
The unitary function $U(.)$ is defined as:

\begin{align}\label{heavyside}
U(\x-\Omega_i)=
\begin{cases}
1,\quad \x \in \Omega_i\\
\\
0, \quad \text{otherwise}
\end{cases}
\end{align}

Hence, the overall optical attenuation coefficient $\mu$ will be a piece-wise
continuous function of spatial coordinate $\x$, where the $\mu_i$ on each
sub-domain $\Omega_i$ is considered to be uncertain.  

\textcolor{black}{We emphasize here that other parameters of the Pennes bioheat equation (e.g. thermal conductivity, perfusion rate, etc.) can also be uncertain in general. However, for  the sake of simplicity, we only considered optical attenuation coefficient, $\mu(\x,t)$, to be uncertain in this manuscript. Note that proposed approach in the following can be easily extended to spatio-temporally varying cases with slight modifications.}

PRF shift is induced by the 
fluctuations in the temperature field~\cite{hindman1966proton}.
We will assume that the changes in the PRF are observed through noise corrupted
$k-$space measurements of a steady state gradient echo sequence:


\begin{equation}\label{obsmodel}
z\propto \int_\Omega \left(M(\x)  e^{-s(\mu,\x)}\right)e^{-2\pi i\overrightarrow{k}.\x} \; d\x+\nu,\quad  \nu \sim \mathcal{N}(0,\sigma)
\end{equation}
where, $\overrightarrow{k}=(k_x,k_y,k_z)$ and
\begin{equation}
s(\mu,\x)=\frac{T_E}{T_2^*(\x)}+i\left\lbrace2\pi\gamma \alpha B_0T_E\Delta u(\mu,\x)+T_E \Delta \omega_0(\x)\right\rbrace
\end{equation}
Transverse magnetization $M(x)$ is determined by the repetition
time, $T_R$, $T_1$-relaxation, and flip angle, $\theta$: 
\begin{equation}\label{mx}
M(\x) = \frac{M_0 \sin \theta \left( 1- e^{-T_R/T_1}\right)}{\left( 1- \cos \theta\ e^{-T_R/T_1}\right)}
\end{equation}
\textcolor{black}{Note that relaxation parameters like $T_1$ can be a function of temperature in general.}
As \eq{obsmodel} represents, observation data is 
polluted with some normally distributed
white noise $\nu$ (with mean 0, and standard deviation $\sigma$). One should
note that complex measurement data $z$ is a function of $k-$space parameters
$\overrightarrow{k}=(k_x,k_y,k_z)$. Table~\ref{sensor_table} summarizes the parameters involved in \eq{obsmodel} to \eq{mx}. 

\begin{table}[htb!]
\centering
\caption{Physical description of the parameters involved in sensor model}\label{sensor_table}
\begin{tabular}{cc}
Parameter (unit) & Description\\
\hline
\hline
$M_0$ (Tesla) & equilibrium magnetization\\
\hline
$\theta$ (degree) & Flip Angle (FA)\\
\hline
$T_1,\ T_2^*$ (sec) & relaxation parameters\\
\hline
$\Delta \omega_0$ (radians)& off-resonance\\
\hline
$T_E$ (sec) & echo time\\
\hline
$\gamma$ ($\frac{\text{rad}}{\text{sec.Tesla}}$) & gyromagnetic ratio of water
\vspace{0.05in}\\
\hline
$\alpha$ ($\frac{\text{ppm}}{^\circ \text{C}}$)& temperature sensitivity\\
\hline
$B_0$ (Tesla) & strength of the magnetic field\\
\hline
\end{tabular}
\end{table}


To be concise in formulation, one can define the observation model as:
\begin{equation}\label{ymodel}
\begin{split}
z & \propto \mathcal{U}(\mu,\overrightarrow{k})+\nu
\qquad
\nu \sim \mathcal{N}(0,\sigma)
\\
\mathcal{U}(\mu,\overrightarrow{k})& =\int_\Omega \left( M(\x)  e^{-s(\mu,\x)}\right)e^{-2\pi i\overrightarrow{k}.\x} d\x
\end{split}
\end{equation}

Our final goal is to accurately estimate the temperature field over the
spatial domain, given a limited number of measurement data $z$. The key
challenge to achieve this goal is to find the \textit{optimal} set of 
frequencies $(k_x,k_y,k_z)$ such that they provide the \textit{maximum
amount of information} for accurate estimation of uncertain parameter $\mu$
and consequently temperature field $u$. 

\section{Methodology}\label{sec:methodology}
\begin{figure}[htb!]
\begin{center}
\includegraphics[trim = 1cm 3cm 8cm 0.5cm, clip,width=3.5in,height=3.5in]{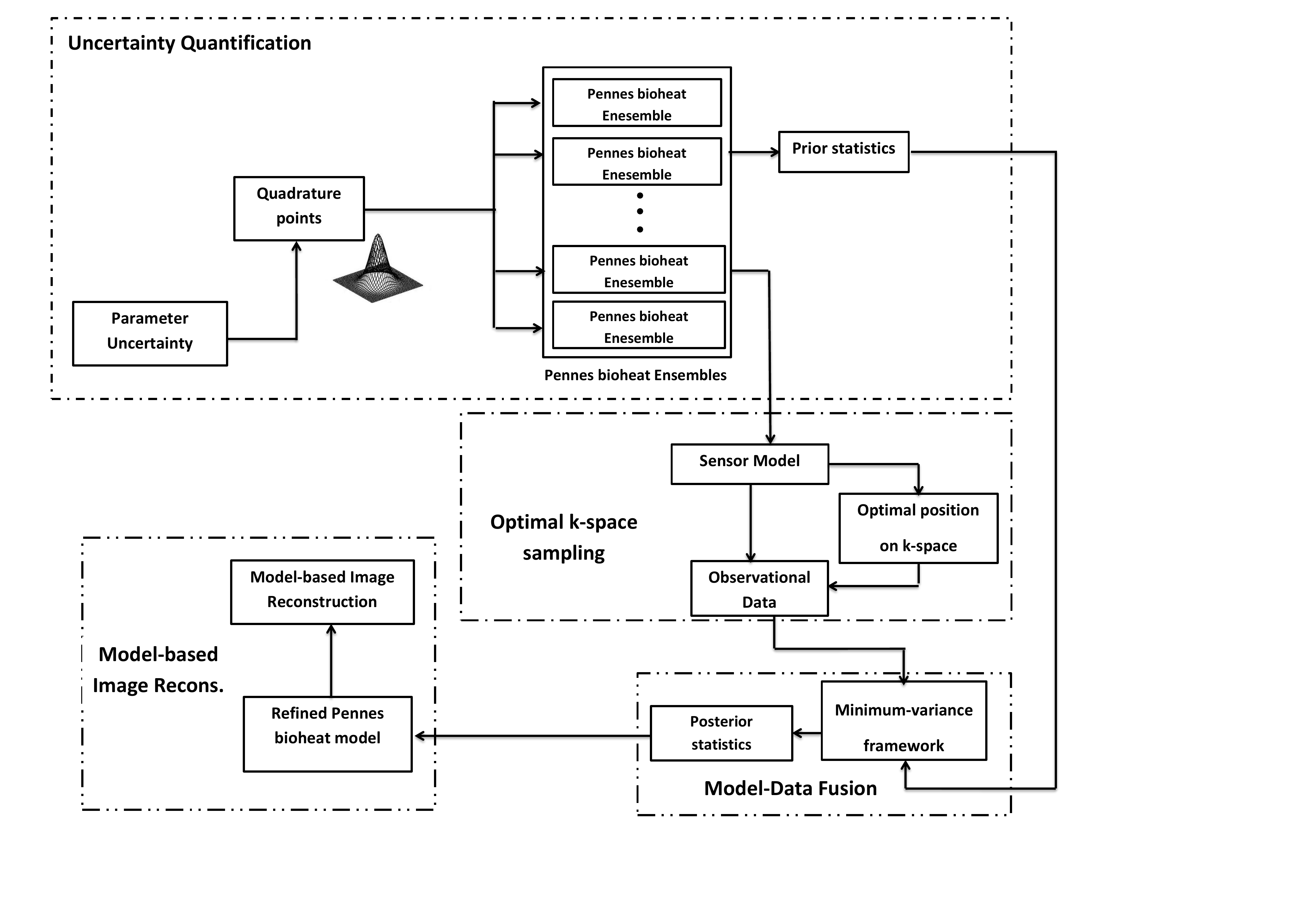}
\caption{Schematic view of estimation process}\label{estimation}\vspace{-0.2in}
\end{center}
\end{figure}

The method presented in this paper consists of four different components
that are combined together to perform the task of the model-based
thermometry image
reconstruction. These components consist of i) Uncertainty Quantification
(UQ), ii) Optimal $k-$space sampling, iii) Model-Data Fusion, and iv)
Model-based Image Reconstruction. Schematic view of the whole estimation
process is shown in \Fig{estimation}. 

As shown in \Fig{estimation}, the process starts with given uncertainties in
the parameters of Pennes bioheat equation. The first step to perform the
estimation is to quantify the effect of these uncertain parameters on
spatio-temporal distribution of the temperature. This is performed by
propagation of a set of quadrature points through the Pennes bioheat
equation. The weighted average of propagated quadrature points is then used to
determine the statistics of the temperature (mean and covariance) over the
spatial domain at a given time~\cite{fahrenholtz2013generalised,madankan_jgcd,madankan_jcp}. 
Hence, utilizing an
\textit{appropriate} and \textit{efficient} set of quadrature points is of
high importance. 
The methodology for
UQ is discussed in further details
in Section \ref{sec_uq}.

Besides precise quantification of uncertainty, having \textit{useful}
measurement data is also highly important to ensure accurate estimation.
This is the focus of this manuscript and is achieved by optimal $k-$space sampling. 
The key goal of optimal
$k-$space sampling is to optimally locate a set of samples over the
$k-$space such that the measurement of \textit{the most useful data} is
ensured. This is achieved by maximizing the mutual information between measurements $z$ and model parameters $\mu$.
Our approximation to the solution of this problem is presented in Section \ref{sec:dyn_data}.

To improve our level of confidence about the present uncertainties in the
Pennes bioheat equation, prior statistics obtained from UQ and observed data
from optimal $k-$space sampling are combined together within a model-data fusion (data
assimilation) framework \cite{madankan_jcp}. 
This will result in \textit{posterior} statistics (e.g.
mean and covariance) of uncertain parameters. Mathematical details of applied data assimilation technique are provided in Section \ref{sec:gpc_minvariance}. 

The obtained posterior statistics are then utilized to \textit{refine} the
parameter values in Pennes bioheat equation, and simulation results provided
by this \textit{updated} mathematical model are used to precisely
reconstruct the spread of temperature over the tissue. 
Detailed description
of this process is provided in Section \ref{sec:reconstruction}.

\textcolor{black}{In here, we should emphasize that the LITT procedure is a dynamic process where the temperature field changes during time. Hence, one would expect the corresponding locations for optimal $k-$space observation to change over the time. However, in this paper we only focus on optimal $k-$space location at a given time step. In fact, temporal change of $k-$space data acquisition points during the LITT procedure is an ongoing subject of research in our group which its result will be published in a separate paper.} 

\section{Uncertainty Quantification}\label{sec_uq}
Applying efficient and accurate tools to quantify the uncertainty is necessary for a precise parameter estimation. This can be achieved by different tools ranging from simple Monte Carlo (MC) method to generalized Polynomial Chaos (gPC) \cite{xiu2002wiener}\cite{wiener1938homogeneous} and the solution of Fokker-Planck-Kolmogorov Equation (FPKE)\cite{jazwinski1970}. The choice of each method for uncertainty quantification depends on the desired level of accuracy and available computational resources. 


In here, the method of quadrature points is utilized to perform the task of uncertainty quantification. In the method of quadrature points, a set of \textit{intelligently} selected points will be propagated through dynamical model and statistics of the output (e.g. mean and variance) are then determined by the weighted average of the model outputs. The major benefit of the quadrature scheme is its significantly less computational cost comparing to Monte Carlo method. Furthermore, the computational complexity involved in quadrature methods is far less than the computational complexities involved in the solution of FPKE and one can easily make use of parallel computing techniques to expedite the procedure. In this paper, due to mentioned benefits of the method of quadrature points with respect to other approaches like MC, gPC, and FPKE, a quadrature based method is being used for quantification of uncertainty in Pennes bioheat equation subjected to parameteric uncertainty.

The method of quadrature points can be viewed as a MC-like evaluation of system of equations, but with sample points selected by quadrature rules. To explain this in more detail, consider $\mu$ to be a time invariant uncertain parameter which can be considered as a function of a random vector $\xii=[\xi_1,\xi_2,\cdots,\xi_{d}]\in \re^{d}$ defined by a pdf $\textit{p}(\xii)$ over the support $\Omega_{\xi}$.

In this way, expected value (mean) of temperature at a given spatio-temporal point $(\x,t)$ can be written as
\begin{align}\label{expx1}
\hat{u}=\Ex{u}{} = \int_{\xii} u(\mu,\x,t)p(\xii)d\xii\simeq \sum\limits_{q}^Mw_q u(\mu(\xii^q),\x,t)
\end{align}
where, $M$ denotes total number of applied quadrature points and $\mu(\xii^q)$ represents $q^{th}$ quadrature point, generated based on the applied quadrature scheme. Similarly, $j^{th}$ order central moments of temperature at each point $(\x,t)$ can be  formulated as
\begin{align}\label{expx1_ctr}
\Ex{(u-\hat{u})^j}{} = \int_{\xii} (u(\mu,\x,t)-\hat{u})^j\ p(\xii)d\xii\simeq \nonumber\\
\sum\limits_{q}^Mw_q\left(u(\mu(\xii^q),\x,t)-\hat{u}\right)^j, \quad j=2,3,\cdots
\end{align}

One can also find moments of $\mathcal{U}$ in similar way, i.e.
\begin{align}\label{expx1U}
m_1=\hat{\mathcal{U}}=\Ex{\mathcal{U}}{} = \int_{\xii} \mathcal{U}(\mu,\x,t)\ p(\xii)d\xii\simeq \nonumber\\
\sum\limits_{q}^Mw_q\ \mathcal{U}(\mu(\xii^q),\x,t)
\end{align}
and 
\begin{align}\label{expx1_ctrU}
m_j=\Ex{(\mathcal{U}-\hat{\mathcal{U}})^j}{} = \int_{\xii} (\mathcal{U}(\mu,\x,t)-\hat{\mathcal{U}})^j\ p(\xii)d\xii\simeq \nonumber\\
\sum\limits_{q}^Mw_q\left(\mathcal{U}(\mu(\xii^q),\x,t)-\hat{\mathcal{U}}\right)^j, \quad j=2,3,\cdots
\end{align}
Hence, the moments at any time can be approximated as a weighted sum of the outputs of simulation runs initiated at the quadrature points generated from the \textit{prior} uncertain parameter distribution. Note that moments for temperature, given by \eq{expx1} and \eq{expx1_ctr}, span over the spatial domain, while moments of $\mathcal{U}$, obtained by \eq{expx1U} and \eq{expx1_ctrU}, span over $k-$space. 

\textcolor{black}{Different types of quadrature schemes like classical Gaussian quadrature rule can be used for uncertainty quantification. However, for a general $m-$dimensional integral, the tensor product of 1-dimensional Gaussian quadrature points results in an undesirable exponential growth of the number of points. In this manuscript, we have used a recently developed Conjugate Unscented Transform (CUT) \cite{venkat} to overcome this drawback of regular quadrature schemes in higher dimensions. Note that CUT points are equivalent with classical Gaussian quadrature points in one dimension. The proposed CUT points are efficient in terms of accuracy while integrating polynomials and yet just employ a small fraction of the number of points used by the traditional Gaussian quadrature scheme. We have also made use of GPU-accelerated parallel computing techniques to expedite the uncertainty quantification process.}

\textcolor{black}{We emphasize here that the method of quadrature points can be used in presence of multiple uncertain parameters in general. However, we only considered the optical attenuation coefficient $\mu$ to be uncertain in this manuscript. The reason for considering only $\mu$ to be uncertain in this manuscript is that among the parameters involved in Pennes Bioheat equation, the temperature distribution is more sensitive to optical attenuation coefficient \cite{fahrenholtz2013generalised}. }

\section{Optimal $k-$space Sampling}\label{sec:dyn_data}
As mentioned earlier, one of the main challenges in MRI technique is to effectively determine the optimal locations on $k-$space such that they provide the best estimate for the image. In other words, we are looking for the measurement data that provide the most confident estimates of uncertain parameters. This is equivalent with finding the measurement data that result in the most reduction of uncertainty in parameter estimates. Based on Information Theory \cite{cover2012elements}, the reduction of uncertainty in one parameter due to the knowledge of the other variable is known as \textit{Mutual Information}. Hence, a good measurement data is the one that provides us the maximum amount of mutual information. Information theoretic approaches and the concept of mutual information has been used in numerous works for optimal sensor placement and sensor management \cite{bourgault2002information,martinez2006optimal,tharmarasa2007large,williams2007approximate,krause2008,choi2010continuous,julian2012distributed,madankan_jae,Madankan_dydess}. The concept of mutual information is also recently being used for optimal measurement design in functional MRI (fMRI) \cite{yan2014linear}. In the following, we first describe the concept of mutual information. Then the mathematical details for finding the optimal locations on $k-$space are presented.

\subsubsection*{Mutual Information as a Measure of Sensor Performance}
According to information theory, the mutual information \cite{cover2012elements} between the model forecast $\mathcal{U}$ and measurement $z$ can be written as:

\begin{equation}\label{mutual_info_u}
I(\mathcal{U}(\mu,\bar{k});z) = \int_{z}\int_{\mathcal{U}} p(\mathcal{U},z)\ln\left(\frac{p(\mathcal{U},z)}{p(\mathcal{U})p(z)}\right)d\mathcal{U} dz
\end{equation}

Since the uncertain parameter $\mu$ is the only source of uncertainty in model forecast $\mathcal{U}$, one can evaluate the mutual information between the uncertain parameter $\mu$ and measurement $z$, instead of \eq{mutual_info_u}:

\begin{equation}\label{mutual_info}
I(\mu;z) = \int_{z}\int_{\mu} p(\mu,z)\ln\left(\frac{p(\mu,z)}{p(\mu)p(z)}\right)d\mu dz
\end{equation}
Using Bayes' rule, $p(\mu,z)$ can be written as

\[p(\mu,z)=p(\mu|z)p(z)\]
Hence $I(\mu;z)$ will be equal to:
\begin{eqnarray}\label{info}
I(\mu;z) = \int_{z} \underbrace{\int_{\mu} p(\mu|z)ln\left(\frac{p(\mu|z)p(z)}{p(\mu)p(z)}\right)d\mu}_{D_{KL}(p(\mu|z))||p(\mu)} p(z) dz
\end{eqnarray}
or,
\begin{equation}
I(\mu;z)=\mathcal{E}_z\left[D_{KL}\left(p(\mu|z))\ ||\ p(\mu)\right)\right]
\end{equation}
  
where, $D_{KL}\left(p(a)\ ||\ p(b)\right)$ denotes the KL distance between two probability density functions $p(a)$  and $p(b)$. Hence, mutual information can be interpreted as the average Kullback-Leiber distance between the prior pdf $p(\mu)$ and the posterior pdf $p(\mu|z)$. By maximizing the mutual information, one inherently maximizes the difference between the prior and posterior distributions of the parameter $\mu$, thus leading to a better measurement and estimate.

\subsection{Optimal Locations on $k-$space}\label{subsec:kspace}
Note that due to dependence of observation data $z$ on $k-$space parameters, the mutual information $I(\mu;z)$ is also a function of $(k_x,k_y,k_z)$ parameters. Hence, our goal is to look for a set of optimal locations $(k_x^j,k_y^j,k_z^j),\ j=1,2,\cdots,N$ such that they maximize the mutual information. This can be mathematically described as the following optimization problem:

\begin{equation}\label{cost1}
\max_K J(K) = \max_{K} I(\mu;\mathbf{z})
\end{equation}
where, $K=\{(k_x^1,k_y^1,k_z^1),\cdots,(k_x^N,k_y^N,k_z^N)\}$ denotes the coordinates of $N$ observations on $k-$space and $\mathbf{z}=\{z_1,z_2,\cdots,z_N\}$ is the set of all $N$ observations in hand. 

We should emphasize here that obtained $k$-space points are completely independent from \textit{noise realizations}, while \textit{statistics of the noise} will affect the location of $k-$space points. For instance, if the standard deviation of noise increases, it will result in a less confident estimate of parameter $\mu$, and vice versa. In extreme case, when standard deviation of noise goes to infinity, i.e. when SNR approaches to zero, there would not be any useful information in $k-$space measurements and proposed algorithm returns the same prior estimates of $\mu$. On the other hand, when SNR goes to infinity (i.e. standard deviation of noise goes to zero), the proposed algorithm will be much more confident about the parameter estimate $\mu$ and correspondingly temperature estimate. 

The major drawback of above optimization problem is its computational intractability for most of the practical applications. Unfortunately, evaluation and optimizing $I(\mu;\z)$ is computationally intractable for most of the practical applications (where the value of $N$ is usually large). Hence, one needs to simplify the original optimization problem and find the approximate solution for optimal location of each observation.

%

\subsection{A Simpler Alternative: Maximizing the Variance}\label{susec:variance}
A simpler alternative in finding useful $k-$space locations is to maximize the entropy (uncertainty) in model outputs $\mathcal{U}$, instead of maximizing the mutual information $I(\mu;z)$ \cite{krause2008}. Based on information theory, the entropy of $\mathcal{U}$, denoted by $h(\mathcal{U})$ is defined as:
%

\begin{equation}\label{entropy}
h(\mathcal{U})=-\int_{\mathcal{U}} ln\left(p(\mathcal{U})\right)p(\mathcal{U})d\mathcal{U}
\end{equation}
where, $p(\mathcal{U})$ is the probability density function of $\mathcal{U}$. 
We emphasize here that this simplification may have its drawback in not giving the \textit{globally} optimal locations on $k-$space, but as we will show in the following, it will result in great simplification of the problem. 

To proceed with maximizing the entropy $h(\mathcal{U})$, we first note that based on Maximum Entropy Principle \cite{cover2012elements}, probability density function $p(\mathcal{U})$ can be parametrized in terms of its central moments as

\begin{equation}\label{pu_entropy}
p(\mathcal{U})=\lim_{l\rightarrow \infty} \left(e^{\sum\limits_{n=0}^l \lambda_n (\mathcal{U}-\hat{\mathcal{U}})^n}\right),\quad \lambda_n \in \re
\end{equation}

where, $\hat{\mathcal{U}}$ is provided by \eqref{expx1U}. By substituting \eq{pu_entropy} in \eq{entropy}, we will have:
$  $
\begin{align}
h(\mathcal{U})=-\int_{\mathcal{U}} \lim_{l\rightarrow \infty} \left(\sum\limits_{n=0}^\infty \lambda_n (\mathcal{U}-\hat{\mathcal{U}})^n\right) p(\mathcal{U})d(\mathcal{U})\nonumber\\
= \lim_{l \rightarrow \infty} \left(-\lambda_0 - \sum\limits_{n=2}^{l} \lambda_n m_n\right) \label{entropy_moments}
\end{align}
where, $m_n$'s are the central moments of $\mathcal{U}$, defined by \eq{expx1_ctrU}. Hence, entropy of $p(\mathcal{U})$ can be described in terms of its central moments, as shown in \eq{entropy_moments}. Now, by approximating $p(\mathcal{U})$ with its first two moments and truncating the above expansion (i.e. letting $l=2$), we have 

\begin{equation}\label{entropy_var}
h(\mathcal{U})\simeq -\lambda_0-\lambda_2 m_2,\quad \lambda_0,\lambda_2 \in \re
\end{equation}
Therefore, to maximize the entropy one can only maximize the variance, i.e.
\begin{eqnarray}\label{entropy_var_max}
\max_{K} h(\mathcal{U}) \simeq -\lambda_0 -\lambda_2 \max_{K}(m_2), \quad \lambda_0,\lambda_2 \in \re
\end{eqnarray}
\textcolor{black}{subject to:
\begin{equation}\label{sparsity_con}
|K^i-K^j|\geq N_d,\quad \forall\ i,j\in 1,2,\cdots,N,\quad i\neq j 
\end{equation}}
where, $|K^i-K^j|=|(k_x^i,k_y^i,k_z^i)-(k_x^j,k_y^j,k_z^j)|$ denotes the distance between the $i^{th}$ and $j^{th}$ locations on $k-$space and $m_2=Var(\mathcal{U})$ is defined from \eq{expx1_ctrU}. Hence, the locations on $k-$space with the highest value of variance for $\mathcal{U}$ are a good approximation of optimal locations for data observations. 

\textcolor{black}{\eq{sparsity_con} is considered to ensure that every distinct pair of $k-$space observations are at least by a distance $N_d$ apart from each other. This constraint is used in order to compensate possible dependencies that can be introduced due to approximation of the original problem. We assumed $N_d=1$ through this manuscript.}

\textcolor{black}{One should note that different order of truncation (i.e. different values of $l$) can be used to approximate the entropy in \eq{entropy_moments}. Clearly, higher order approximation (greater values of $l$) leads to more accurate approximation of entropy. However, the downside of using higher order terms is that one needs to find the corresponding $\lambda_n$ coefficients for each term. Finding corresponding $\lambda_n$ coefficients requires solving an optimization problem which could increase computational cost of the whole procedure. Hence, lower order approximation of entropy is of more interest due to real time applications of the proposed technique.}

%
The intuition behind the idea of maximizing the variance is that the points on $k-$space with lower value of variance are less sensitive to model perturbations (resulting from parameter uncertainties) and vice versa. Hence, it is better to select the points on $k-$space with highest sensitivity with respect to model uncertainties. For instance, a point with zero variance on $k-$space will not be a good candidate for data observation since no matter what the values of uncertain parameters are, model output will always be the same at that specific point. On the other hand, a point on $k-$space with large value of variance means that the model output at that location is \textit{highly sensitive} to model uncertainties. Hence, a measurement at that $k-$space location would be of more interest. Therefore, the points with highest values of variance of $\mathcal{U}$ are \textit{better candidates} for data observations. 


\subsection{Adaptation to Practical $k-$space Sampling}\label{susec:var_kspace}
It is well known that practically, $k-$space sampling occurs along readout lines, rather than just sparse points on $k-$space. Hence, one needs to modify the idea of variance maximization in order to make it applicable to practical $k-$space sampling. To do so, we propose to perform the $k-$space sampling along the readout lines which their points posses the higher value of variance. For instance, in 2 dimensions, we have:

\begin{equation}\label{var2d}
\max_{K_y} \sum_{k_x=1}^{N_x} m_2(k_x,k_y)
\end{equation}
where, $K_y=\{k_y^1,\cdots,k_y^N\}$ and $N_x$ denotes the total number of points along the $k_x$ axis on $k-$space. Similarly, in $3-$dimensional $k-$space sampling, one can write
\begin{equation}\label{var3d}
\max_{K_{xy}} \sum_{k_z=1}^{N_z} m_2(k_x,k_y,k_z)
\end{equation}
where, $K_{xy}=\{(k_x^1,k_y^1),\cdots,(k_x^N,k_y^N)\}$ and $N_z$ denotes the total number of points along the $k_z$ axis on $k-$space. 

The major benefit of using \eq{var2d} and \eq{var3d} instead of maximizing \eq{entropy_var_max} is that it provides a more practical scheme for $k-$space undersampling by providing the \textit{most informative readout lines}, rather than the most informative points. In addition, less computational effort is required to solve  \eq{var2d} and \eq{var3d} comparing to \eq{entropy_var_max}. Note that after finding the readout lines through maximizing \eq{var2d} and \eq{var3d} all the points along those readout lines are used for model - data fusion.

\section{Model-Data Fusion}\label{sec:gpc_minvariance}
\textcolor{black}{Once, the most useful measurement data are extracted through the proposed technique, one can combine these data with model predictions to improve our level of confidence about the uncertain parameters and consequently temperature estimates.}

Various approaches exist to perform the model-data fusion. Here, we employ a
linear unbiased minimum variance estimation method for this purpose. In
detail, the minimum variance technique is used to minimize the trace of the
posterior parameter covariance matrix: \begin{equation}\label{costfun}
J = \min\limits_{\mu}\ Tr\left[\Ex{(\mu-\Ex{\mu}{})(\mu-\Ex{\mu}{})^T}{}\right]
\end{equation}
where, $\mu=[\mu_1,\mu_2,\cdots,\mu_d]^T$ denotes optical attenuation
coefficients for $d$ tissues. Note that the minimum variance formulation is
valid for any pdf, although the formulation makes use of only the mean and
covariance information. It provides the maximum \textit{a posteriori} estimate when
model dynamics and measurement model are linear and state uncertainty is
Gaussian. Minimizing the cost function $J$ subject to the constraint of
being an unbiased estimate,  and using linear updating, results in the first
two  moments of the posterior distribution
\cite{gelb,madankan_jcp,madankan_jgcd}:
\begin{align}
\hat{\mu}^+&=\hat{\mu}^- + \mathbf{K}[{\z}-\Ex{\mathcal{U}(u,\mu)}{-}]\label{minvar_mean}\\
\boldsymbol{\Sigma}^+& = \boldsymbol{\Sigma}^- - \mathbf{K} \boldsymbol{\Sigma}_{\mathcal{U}\mathcal{U}}\mathbf{K}^T\label{minvar_var}
\end{align}
where, $\z=[z_1,z_2, \cdots, z_N]^T$ is the measurement vector of $N$
observations and the gain matrix $\mathbf{K}$ is given by \begin{equation}
\mathbf{K} = \boldsymbol{\Sigma}_{\mu z}\left(\boldsymbol{\Sigma}_{\mathcal{U}\mathcal{U}}^-+\mathbf{R}\right)^{-1}\label{Kgain}
\end{equation}
Here, 
$\hat{\mu}^-$ and $\hat{\mu}^+$ represent prior and posterior value of the
mean for the parameter vector $\mu$, respectively:
\begin{eqnarray}\label{pc_prior}
{\hat{\mu}}^-&\triangleq\Ex{\mu}{-} = \int_{\xii}\mu^-(\xii) p(\xii)d\xii\\
\label{gpc_postmean}
\hat{\mu}^+&\triangleq\Ex{\mu}{+} = \int_{\xii}\mu^+(\xii) p(\xii)d\xii
\end{eqnarray}
Similarly, the prior and posterior covariance matrices
$\boldsymbol{\Sigma}^-$ and $\boldsymbol{\Sigma}^+$ can be written as:
\begin{align}\label{pc_prior}
{\boldsymbol{\Sigma}}^-\triangleq\Ex{(\mu-\hat{\mu}^-)(\mu-\hat{\mu}^-)^T}{-} \\
\label{gpc_postvar}
{\boldsymbol{\Sigma}}^+\triangleq\Ex{(\mu-\hat{\mu}^+)(\mu-\hat{\mu}^+)^T}{+} 
\end{align}
Also, diagonal matrix $\mathbf{R}$ denotes the measurement error covariance
matrix in \eq{Kgain} which encapsulates the measurement's inaccuracies. Note
than $\mathbf{R}_{ii}=\sigma^2$.
The matrices $\boldsymbol{\Sigma}_{\mu z}$ and $\boldsymbol{\Sigma}_{\mathcal{U}\mathcal{U}}$ are defined as:
\begin{align}
\hat{\mathcal{U}}^-\triangleq \Ex{\mathcal{U}(u,\mu)}{-}=
\int_{\xii} \underbrace{\mathcal{U}(u^-(\xii),\mu^-(\xii))}_{\mathcal{U}} p(\xii)d\xii\label{h_k}\\
\boldsymbol{\Sigma}_{\mu z}  \triangleq\Ex{(\mu-\hat{\mu})(\mathcal{U}-\hat{\mathcal{U}}^-)^T}{-} 
\label{Pzy}\\
\boldsymbol{\Sigma}_{\mathcal{U}\mathcal{U}}^- \triangleq \Ex{(\mathcal{U}-\hat{\mathcal{U}}^-)(\mathcal{U}-\hat{\mathcal{U}}^-)^T}{-} 
\label{Pzz}
\end{align}
where, the expectation integrals in \eq{h_k}, \eq{Pzy}, and \eq{Pzz}
can be computed by  suitable quadrature rules.

\section{Model-based Image Reconstruction}\label{sec:reconstruction}
After accurate estimation of uncertain parameters using the proposed
approach, one can simply \textit{refine} the Pennes bioheat model using
obtained posterior mean of uncertain $\mu$, \textit{i.e.}

\begin{align}
\raggedright
\rho c \left(\frac{\partial u}{\partial t}+\nabla u. \mathbf{v}\right)-\nabla.({\lambda(u,\x)}\nabla u)+&\nonumber
\\{\omega(u,\x)}c_{blood}(u-u_a)=Q_{laser}(\x,t), \quad \forall \x \in \Omega &\nonumber\\
Q_{laser}(\x,t)=P(t)\hat{\mu}^{+^2}\frac{e^{-\hat{\mu}^+||\x-\x_0||}}{4\pi ||\x-\x_0||}\nonumber\\
u(\x,0)=u_0(\x),\quad \x \in \Omega & \nonumber\\
u(\x,t)=u_D(\x),\quad \x \in \partial \Omega_D & \label{pennes_post}\\
-\lambda(u,\x)\nabla.u(\x,t)=g_N(\x),\quad \x \in \partial \Omega_N & \nonumber\\
-\lambda(u,\x)\nabla.u(\x,t)=h(u-u_{\infty}),\quad \x \in \partial \Omega_c & \nonumber
\end{align}
where, $\hat{\mu}^+$ is the posterior mean of uncertain parameter $\mu$,
obtained from the model-data fusion. Generated temperature field from
\eq{pennes_post}  is an accurate estimate of the reconstructed image from
$k-$space.  Note that one can also use the posterior variance of $\mu$ to
find the variance in generated temperature field by using \eq{pennes}, \eq{expx1_ctr}, and \eq{minvar_var}.

\section{Results}\label{sec:num_sim}
To validate performance of the proposed approach, an \textit{in silico} and two \textit{in vivo} experiments are considered. 
\textit{In silico} experiment demonstrates performance of the proposed approach in simulation of temperature dispersion
in human brain over a $3-$dimensional space. While the \textit{in vivo}
experiments describe planar temperature image reconstruction
in human brain and 3-dimensional temperature imaging of agar phantom. A detailed description of each of
these examples is provided in the following.

\subsection{Volumetric Image Reconstruction}
We initially study the performance of the proposed technique in
a mathematical phantom.
{ \color{black}
Fully sampled volumetric temperature imaging of MRgLITT is simulated. 
A Gaussian mixture model is assumed to perform the tissue segmentation~\cite{avants2011reproducible}. K-means clustering is used to initialize a 4 tissue intensity model representing grey matter, white matter, cerebral spinal fluid, and enhancing lesion on contrast enhanced $T_1$-weighted imaging. The segmentation problem is solved using expectation maximization of a parametric formulation of the mixture model. Piecewise homogeneity across tissue types is enforced using Markov random fields. The segmentation is registered to the thermal imaging using DICOM coordinate and orientation information. } \Fig{brain} shows the position of the tumor in brain for this
example.  
\begin{figure}[htb!!!]
\centering
\vspace{-0.1in}
\begin{tabular}{c}
\hspace{0.1in}\subfigure[]{\includegraphics[width=2.6in,trim = 0in -0.6in 0in 0in,clip]{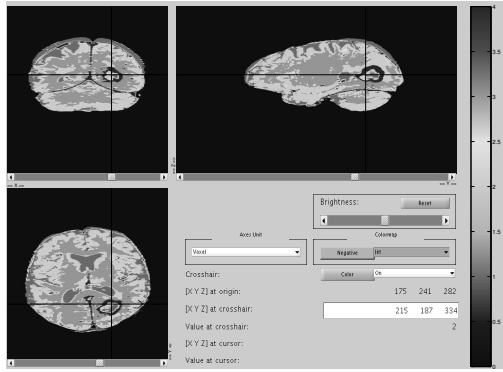}\label{brain}}
\vspace{-0.15in}\\
\subfigure[]{\includegraphics[width=3.3in,height=2.25in,trim = 0in 0.1in 0in 0in,clip]{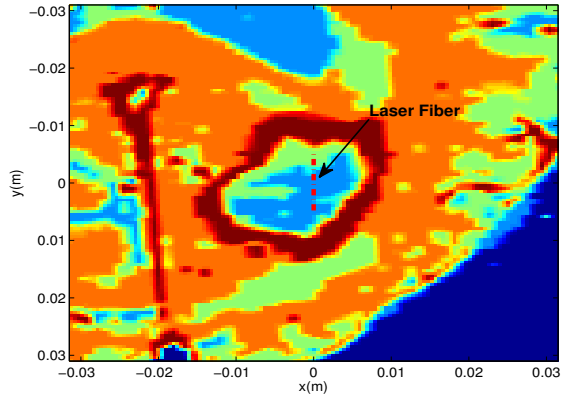}\label{tissue3D}}
\end{tabular}
\vspace{-0.1in}
\caption{Volumetric Image Reconstruction: a) Position of the 
tumor in the brain, projected on different planes, b) Segmentation of the tissue over 
ROI, as viewed over transverse plane (at $36^{th}$ slice in $z$ direction). Four tissue types are shown in
different colors: CSF: light blue, grey matter: green, white matter: orange,
tumor: red. Dashed red lines shows position of laser fiber.}\label{brain} 
\vspace{-0.1in}
\end{figure}

\Fig{tissue3D} illustrates the tissue segmentation in one of the slices within
the Region Of Interest (ROI). As \Fig{tissue3D} shows, there exist four tissue types within the ROI: 
Cerebrospinal Fluid (CSF), grey matter, white matter, and tumor tissue.
We consider the corresponding optical attenuation for each of these
tissues to be uncertain. Hence, we have four uncertain parameters. The
optical attenuation coefficient for the CSF, denoted by $\mu_1$, is assumed
to be uniformly distributed between $10\ \frac{\text{1}}{\text{m}}$ and $300\
\frac{\text{1}}{\text{m}}$, \textit{i.e.}, $\mu_1 \sim \mathcal{U}[10,300]$, while all the other optical
attenuation coefficients are assumed to be uniformly distributed between
$10\ \frac{\text{1}}{\text{m}}$ and $400\ \frac{\text{1}}{\text{m}}$,\textit{ i.e.}, $\mu_i \sim \mathcal{U}[10,400],\
i=2,3,4$. All the other physical properties like  conductivity,
perfusion, etc,  are assumed to be known and are given in Table \ref{pennes_para}. For the LITT process, we assumed the laser power to be 11.5 watts and heating time period is considered to be $90$ sec.

\begin{table}[htb!!]
\caption{
\color{black}
Constitutive
Data used in numerical simulations ~\cite{welch1984thermal,duck1990}}\label{pennes_para}
\centering
\begin{tabular}{ccccc} \hline
$\Lambda $ $ \frac{W}{ m \cdot K}$ & $\omega$ $\frac{kg}{m^3 s}$ &  $\rho$
$\frac{kg}{m^3}$ &   $c_{blood}$ $ \frac{J}{kg \cdot K}$ &  $c$
$\frac{J}{kg \cdot K}$ \\ \hline\hline
          0.527              &             9.0           &  1045
&            3840                      &                  3600          \\
\hline
\end{tabular}
\end{table}

{\color{black}
To perform the estimation process, we have generated fully sampled synthetic measurement
data using a \textit{random realization} of the optical attenuation coefficient
$\mu$. In detail, we used
$(\mu_1,\mu_2,\mu_3,\mu_4)=(111.39,218.75,383.01,385.96)\;
\frac{1}{\text{m}}$ 
to generate the temperature field through the Pennes bioheat model.
The obtained temperature
field was then used in \eq{obsmodel} to generate the corresponding
\textit{fully-sampled} measurement data. The Signal to Noise Ratio (SNR) is assumed to be 50 in this experiment.
} Table \ref{sensor_para} shows the numerical values of the parameters
involved in the sensor structure.
\begin{figure}[htb!]
\centering
\vspace{-0.1in}
\begin{tabular}{c}
\subfigure[]{\includegraphics[width=3.4in,trim = 0in 5in 0in 0.63in,clip]{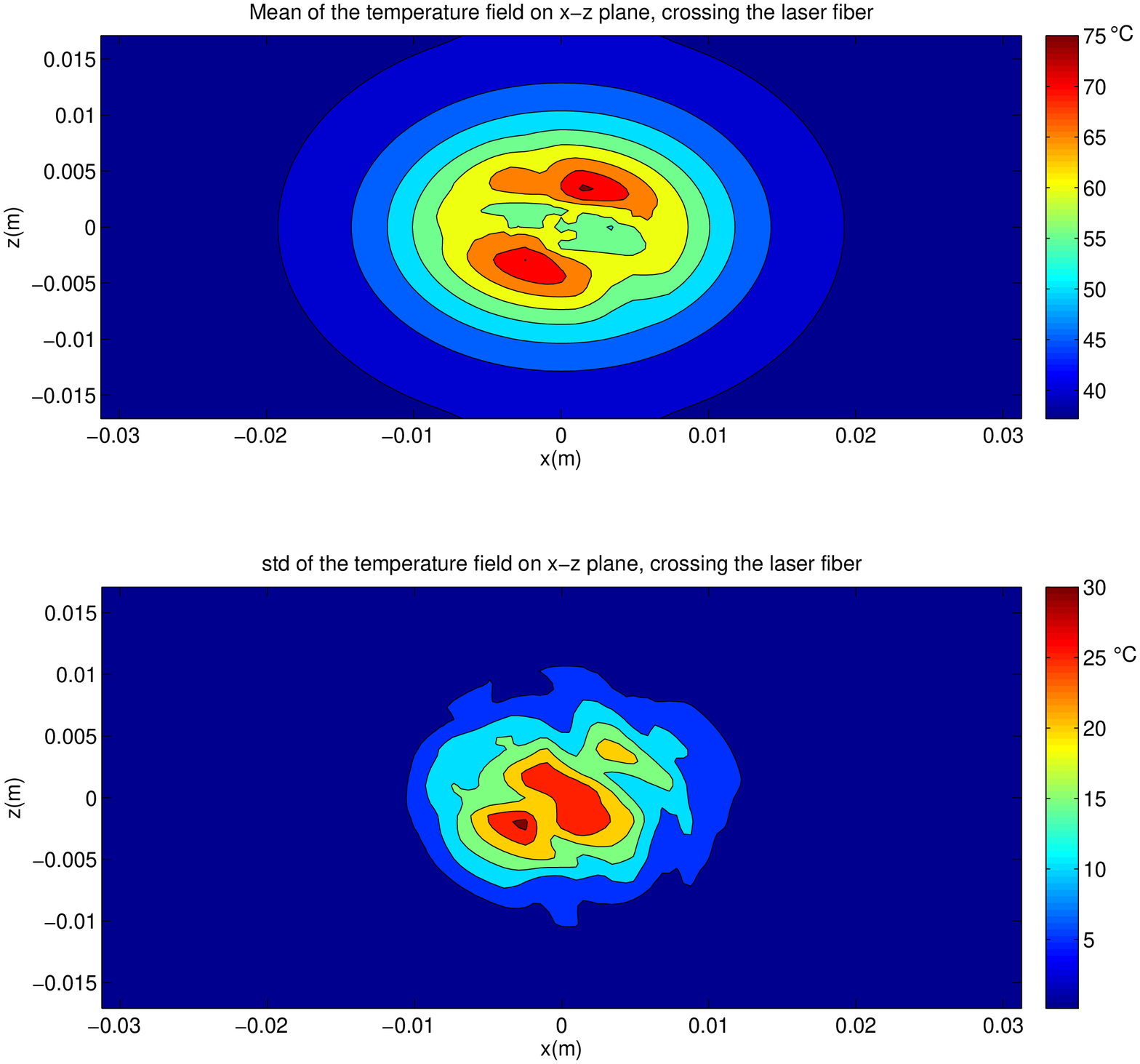}} \\
\hspace{0in}\subfigure[]{\includegraphics[width=3.4in,trim = 0in 0in 0in 5.6in,clip]{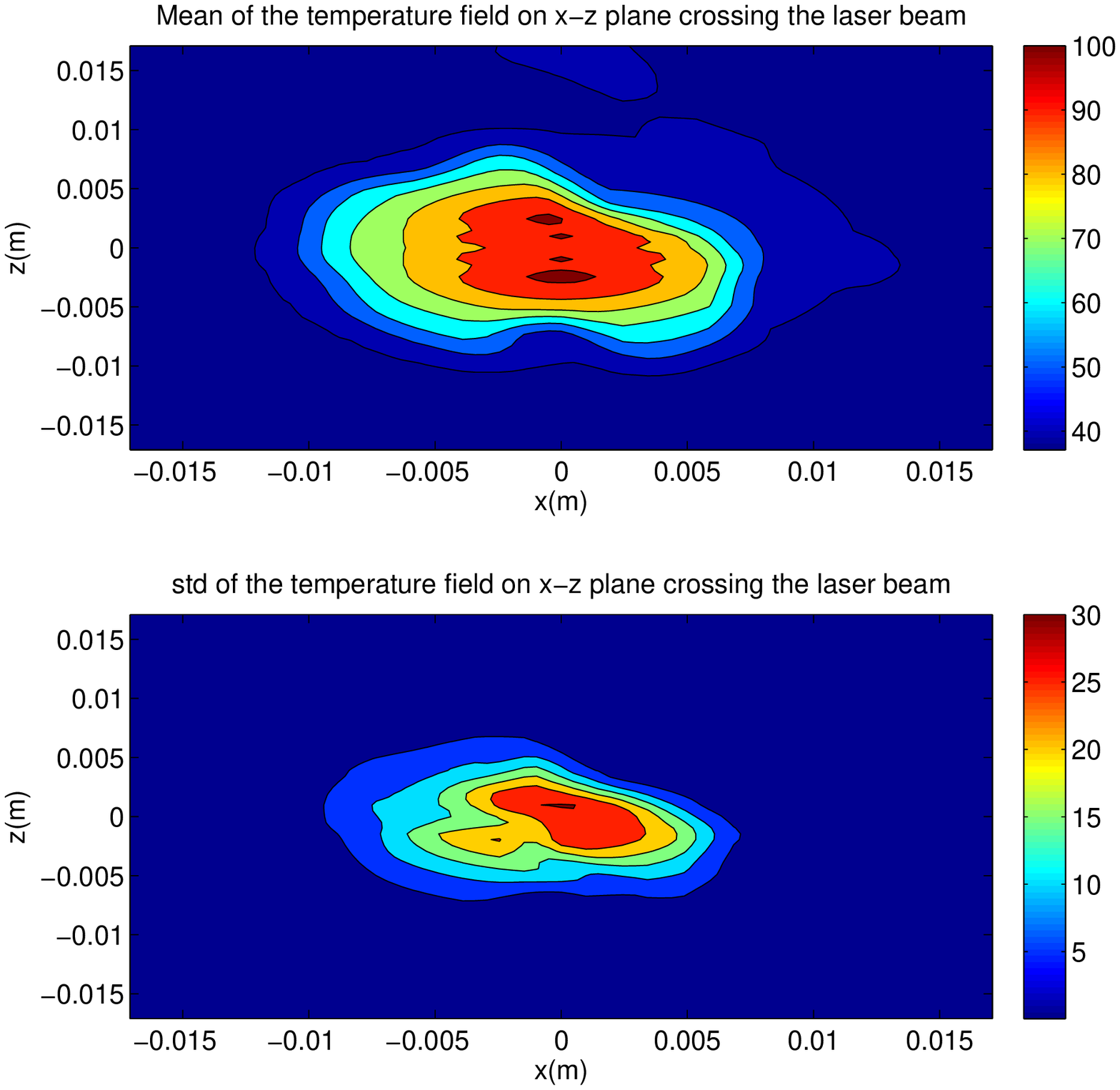}} 
\end{tabular}
\caption{Volumetric Image Reconstruction: a) mean ($^\circ$C) and b) standard deviation of temperature ($^\circ$C) over the ROI crossing through the laser fiber, projected \textit{x-z} plane.}\label{pennesUQ}
\vspace{-0.15in}
\end{figure}

\begin{table}[htb!]
\centering
\caption{Volumetric Image Reconstruction: parameter values involved in sensor model}\label{sensor_para}
\begin{tabular}{cc}
Parameter & Value\\
\hline
\hline
$\theta$ (rad) & $\pi/3$ \\
\hline
$T_1$ (sec.) & 4.31 (CSF), 1.035 (grey), 0.63 (white), 0.8 (tumor) \\
\hline
$T_2^*$ (ms) & 10 (CSF), 70 (grey), 100 (white), 80 (tumor)\\
\hline
$T_R$ (sec.) & 0.544 \\
\hline
$\Delta \omega_0$ (rad) & 0\\
\hline
$T_E$ (ms) & $25$ \\
\hline
$\gamma$ (MHz/T) & 42.58 \\
\hline
$\alpha$ (ppm/C) & $-0.0102$ \\
\hline
$B_0$ (T) & $1.5$ \\
\hline
\end{tabular}
\end{table}

A set of 161 Conjugate Unscented Transform (CUT) points \cite{venkat} (which accurately evaluate upto $8^{th}$ order polynomial integrals) are used to capture the uncertainty in the temperature field due to the presence
of uncertain optical attenuation coefficients. \Fig{pennesUQ} illustrates
the mean and standard deviation of the temperature field over an $x-z$
slice co-planar with the laser fiber.  

\begin{figure*}[htb!!!]
\vspace{-0.1in}
\centering
\subfigure[]{\includegraphics[width=3in]{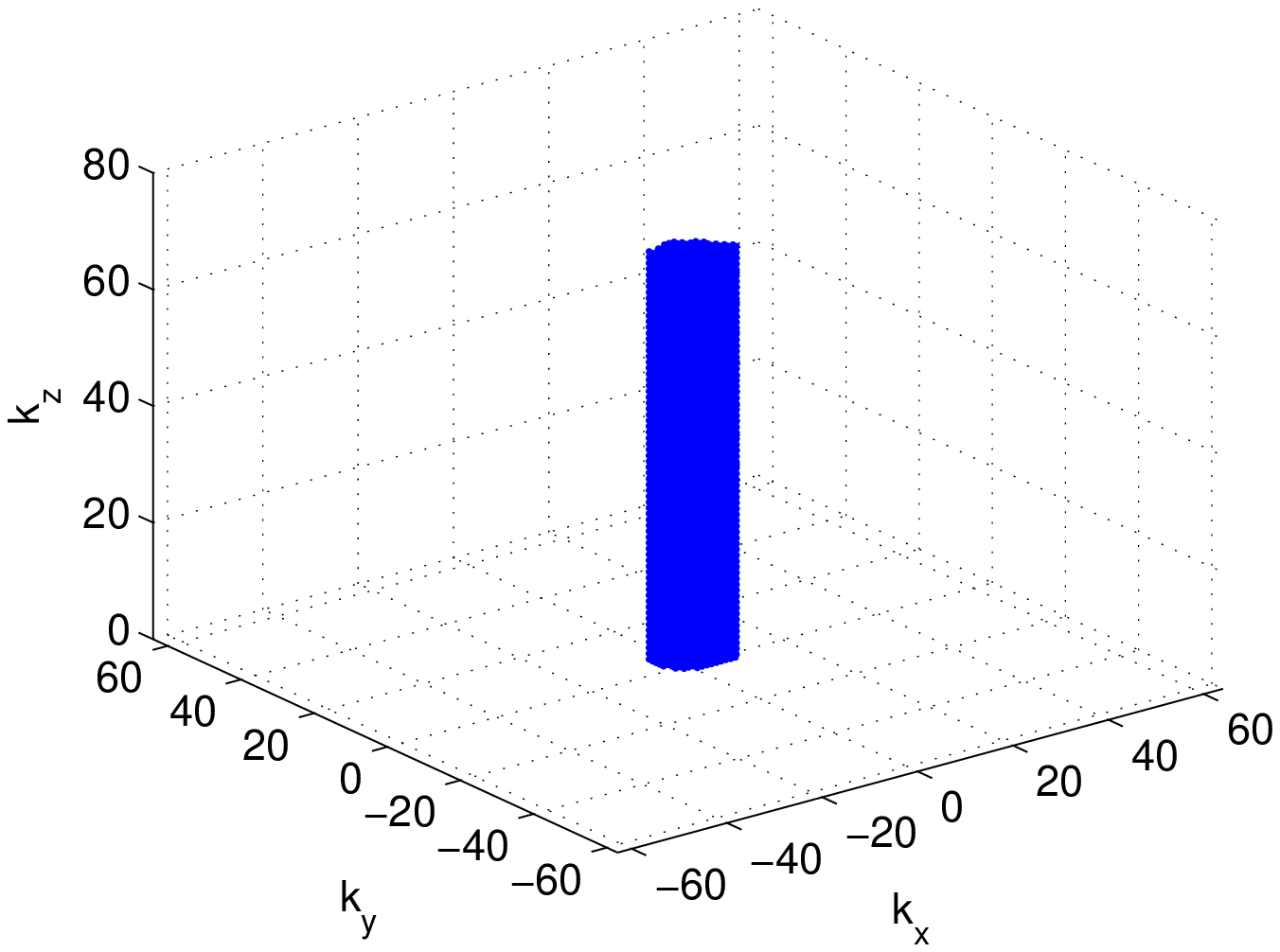}\label{kpoints3D}}
\hspace{0.3in}
\subfigure[]{\includegraphics[width=3in]{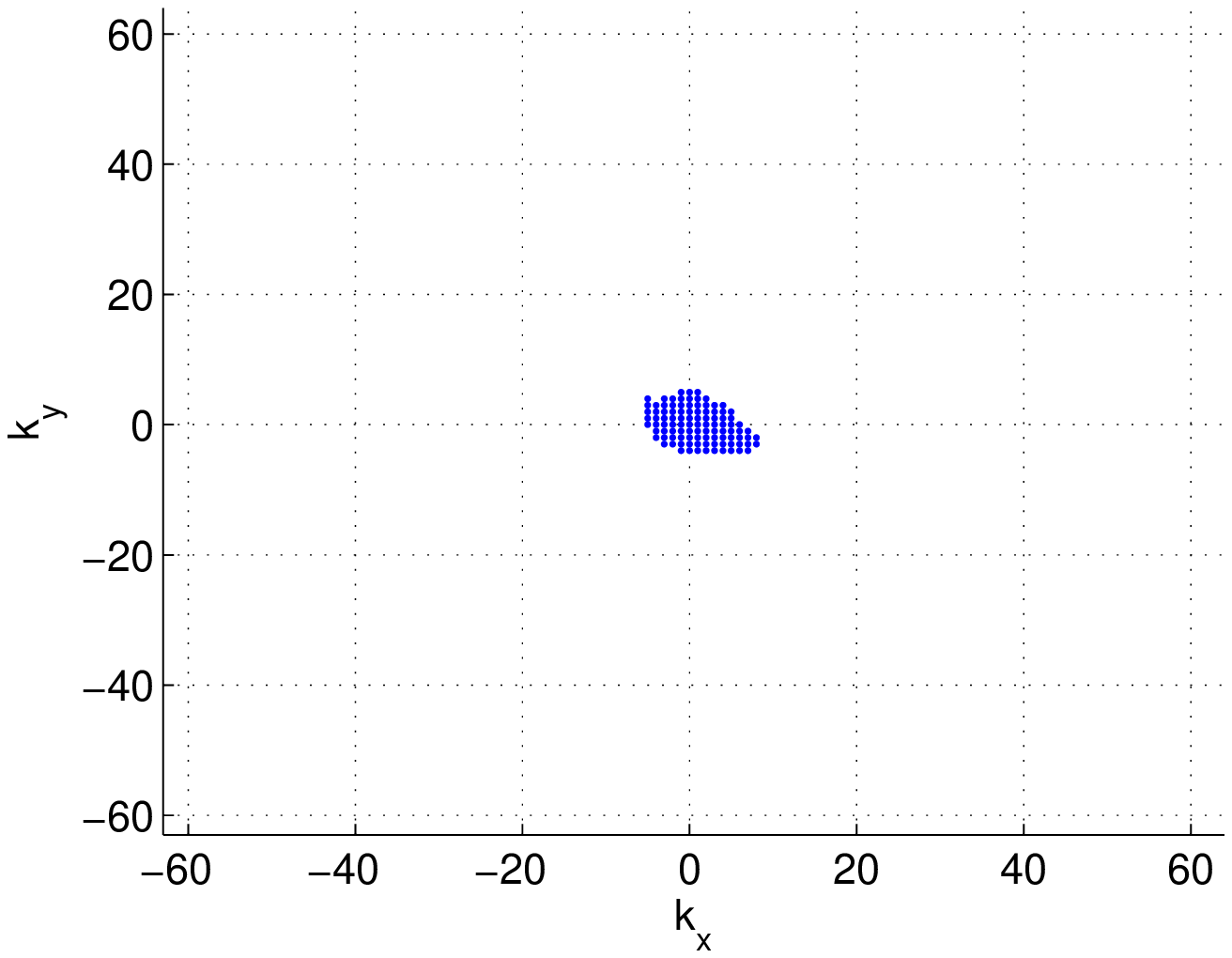}\label{kpointsxy}}
\vspace{-0.1in}
\caption{Volumetric Image Reconstruction: Applied points for data observation a) 3-dimensional view b) projection of the points on $k_x-k_y$ plane. Readout lines are the lines passing through these 100 points.}\label{kpoints}
\vspace{-0.2in}
\end{figure*}

As discussed in Section \ref{sec:dyn_data}, the variance maximization technique
is utilized to approximate  $k-$space location with the highest information
content for the data acquisition. 
{\color{black}
\Fig{kpoints} illustrates 100 lines in $k-$space with the highest
variance.}
Most of the points lie in center of
the $k-$space. 
The number of points projected on the $k_x-k_y$ plane is
very small comparing to total number of points used in data observation.
For instance, in case of using 100 lines on $k-$space for data
observation, the proposed technique makes use of only $\frac{100}{128^2}\simeq 0.6\%$ of the points on $k_x - k_y$ plane. Note that this is equivalent with almost 164 times speed up in comparison with full-sampling of $3-$dimensional $k-$space!
This can be seen in \Fig{kpointsxy}.

\begin{figure}[htb!!!]
\centering
\vspace{-0.1in}
\includegraphics[width=3.4in]{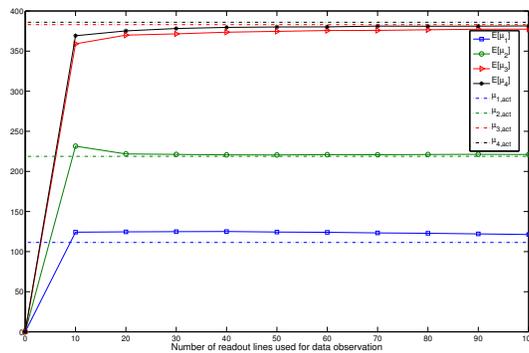}
\vspace{-0.15in}
\caption{
Volumetric Image Reconstruction: Convergence of parameter estimates
versus different number of $k-$space readout lines, dashed lines represent the actual value of corresponding
attenuation coefficient.}\label{mu_cvr} 
\vspace{-0.1in}
\end{figure}

{\color{black}
We performed the model-data fusion using acquired $k-$space
samples with the predicted highest variance content. 
Within the context of the minimum variance framework, Section
\ref{sec:gpc_minvariance},
posterior statistics of optical
attenuation coefficients were then found by merging these selected measurements ($k-$space
samples) with model predictions.} 
\Fig{mu_cvr} illustrates the convergence of the posterior
estimate of optical attenuation coefficients versus different 
numbers of measurement data, ranging from 0 to 100 lines.
As it is shown in \Fig{mu_cvr}, posterior estimates of parameters,
converge to their actual value
(denoted by dashed lines) by increasing the number of data
observations.

\begin{figure}[htb!!!]
\centering
\includegraphics[width=3.2in,height=2.25in]{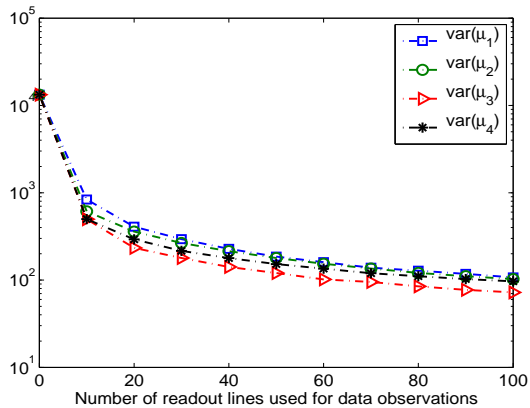}
\vspace{-0.1in}
\caption{Volumetric Image Reconstruction: Variance of parameter estimates
versus different number of $k-$space samples.}\label{mu_var_cvr}
\vspace{-0.1in}
\end{figure}

\Fig{mu_var_cvr} shows the variance for each of the parameter estimates versus
different number of readout lines. As expected, the variance reduces with increasing the number of observations. Note that
the associated variance of posterior parameter estimates is much smaller than
prior variance of the parameters, \textit{i.e.} applied model-data fusion have
significantly increased our level of confidence about the optical
attenuation parameters.  

\begin{figure}[htb!]
\centering
\hspace{0in}{{\includegraphics[width=3.4in]{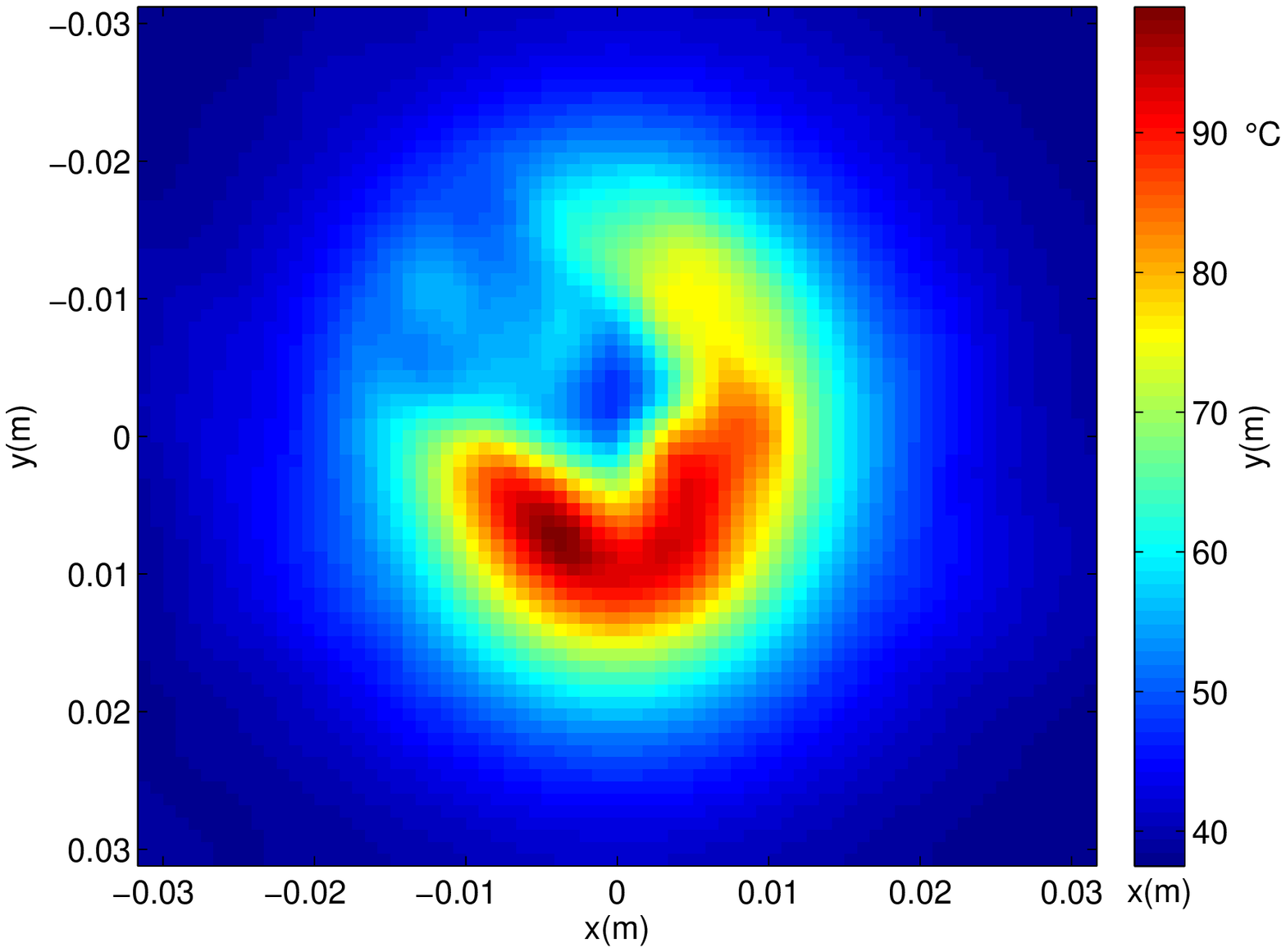}}\label{temp_est}} \\
\caption{Volumetric Image Reconstruction: Posterior estimate of
temperature field ($^\circ$C), using only 100
lines through $k_z$, over the
$40^{th}$ slice in $z$ direction 
}\label{temp_fig} 
\end{figure}

Posterior estimates of $\mu$ are then used in the Pennes bioheat equation to estimate the
true temperature field. {\color{black} \Fig{temp_fig} illustrates temperature field estimate in a representative  $x-y$ plane.}
%

\subsection*{Comparison with Other Under-Sampling Schemes}
\begin{figure*}[htb!!!]
\vspace{-0.1in}
\centering
\begin{tabular}{cc}
\subfigure[]{\includegraphics[width=3in]{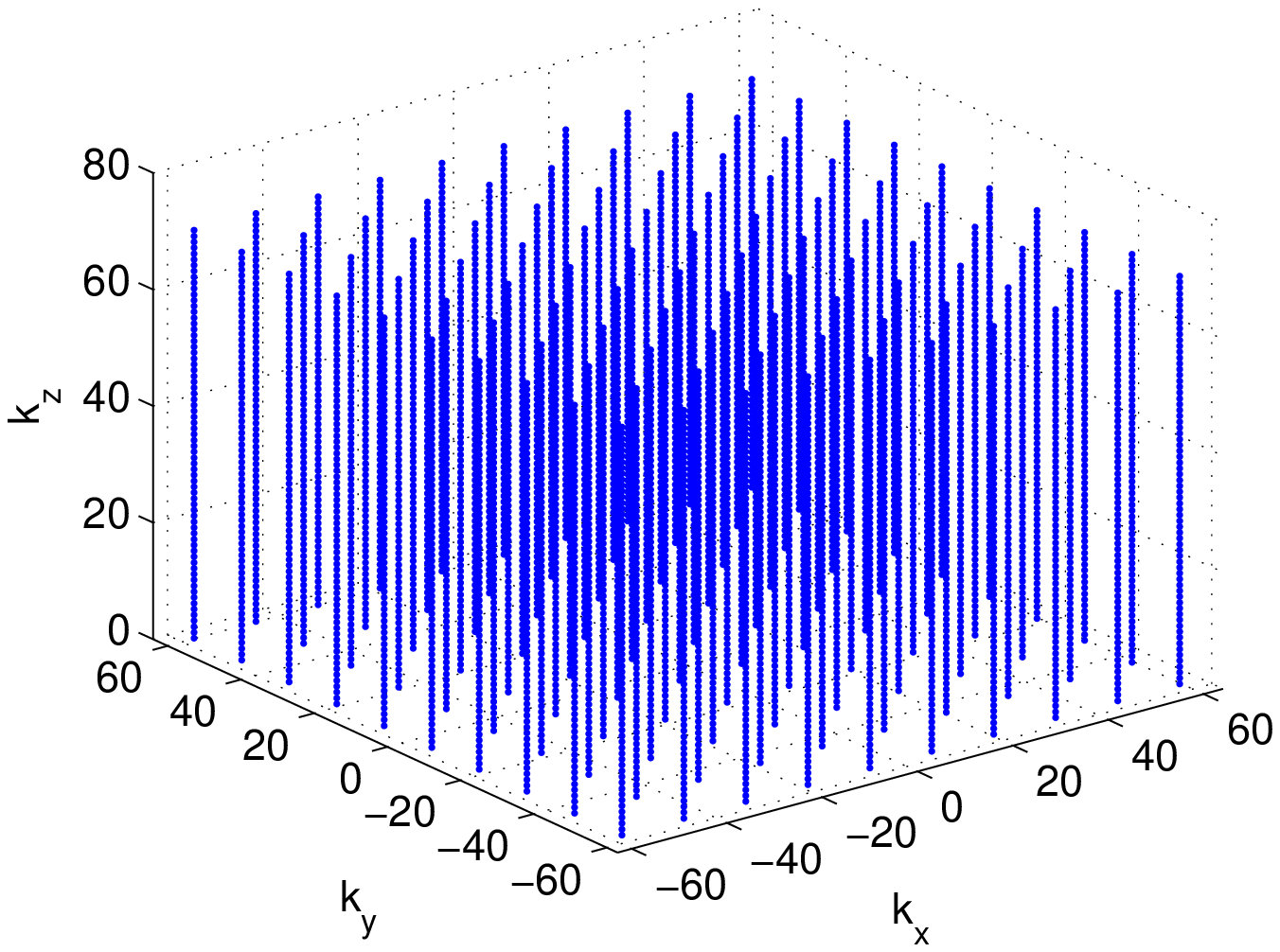}\label{kpoints3Dp}}&
\subfigure[]{\includegraphics[width=3in]{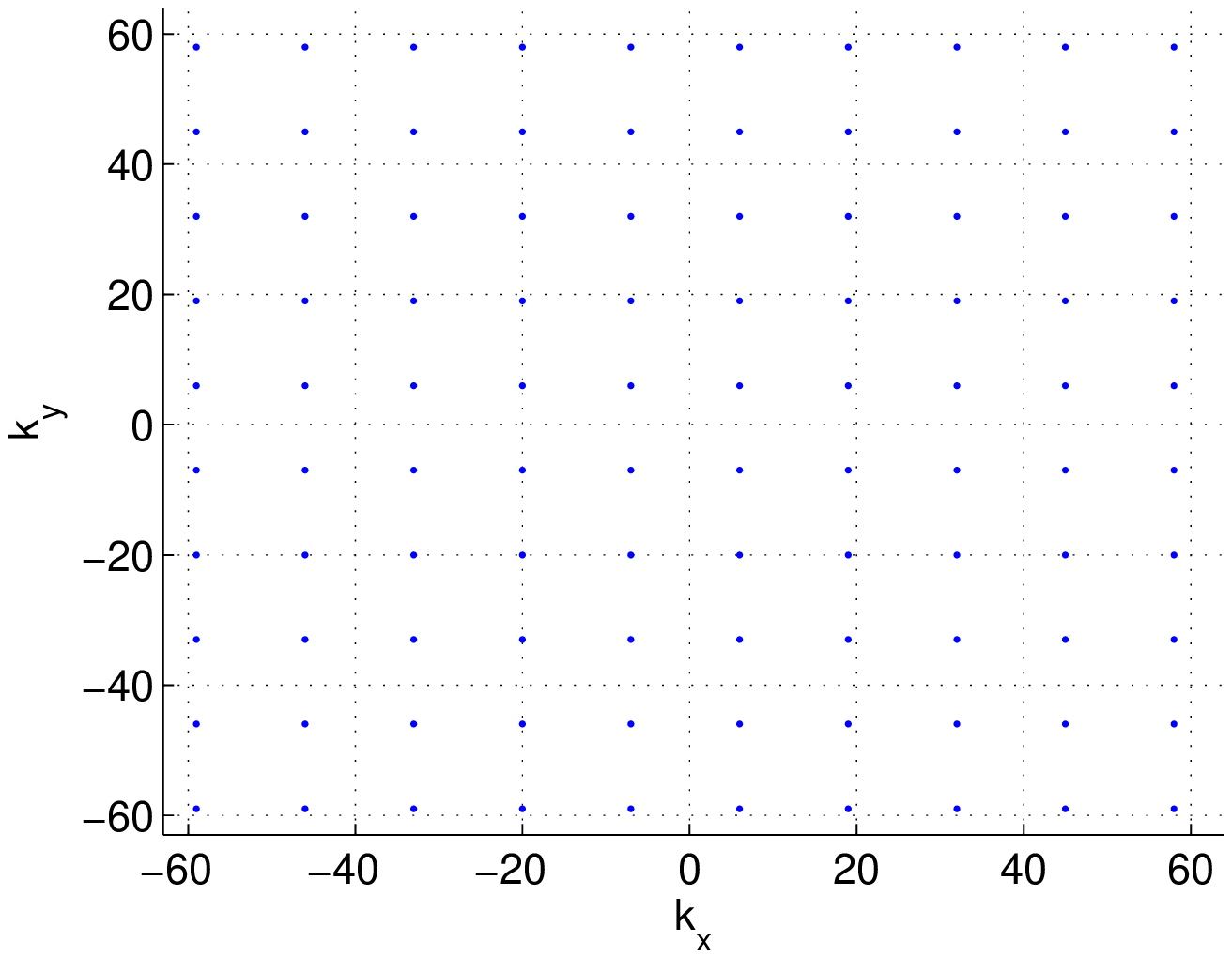}\label{kpointsxyp}}
\end{tabular}
\vspace{-0.1in}
\caption{Volumetric Image Reconstruction: Acquired points for data observation (based on rectilinear undersampling) a) 3-dimensional view b) projection of the points on $k_x-k_y$ plane. Readout lines are the lines passing through these 100 points.}\label{rect3d}
\vspace{-0.1in}
\end{figure*}

\begin{figure*}[htb!!!]
\vspace{-0.1in}
\centering
\begin{tabular}{cc}
\subfigure[]{\includegraphics[width=3in]{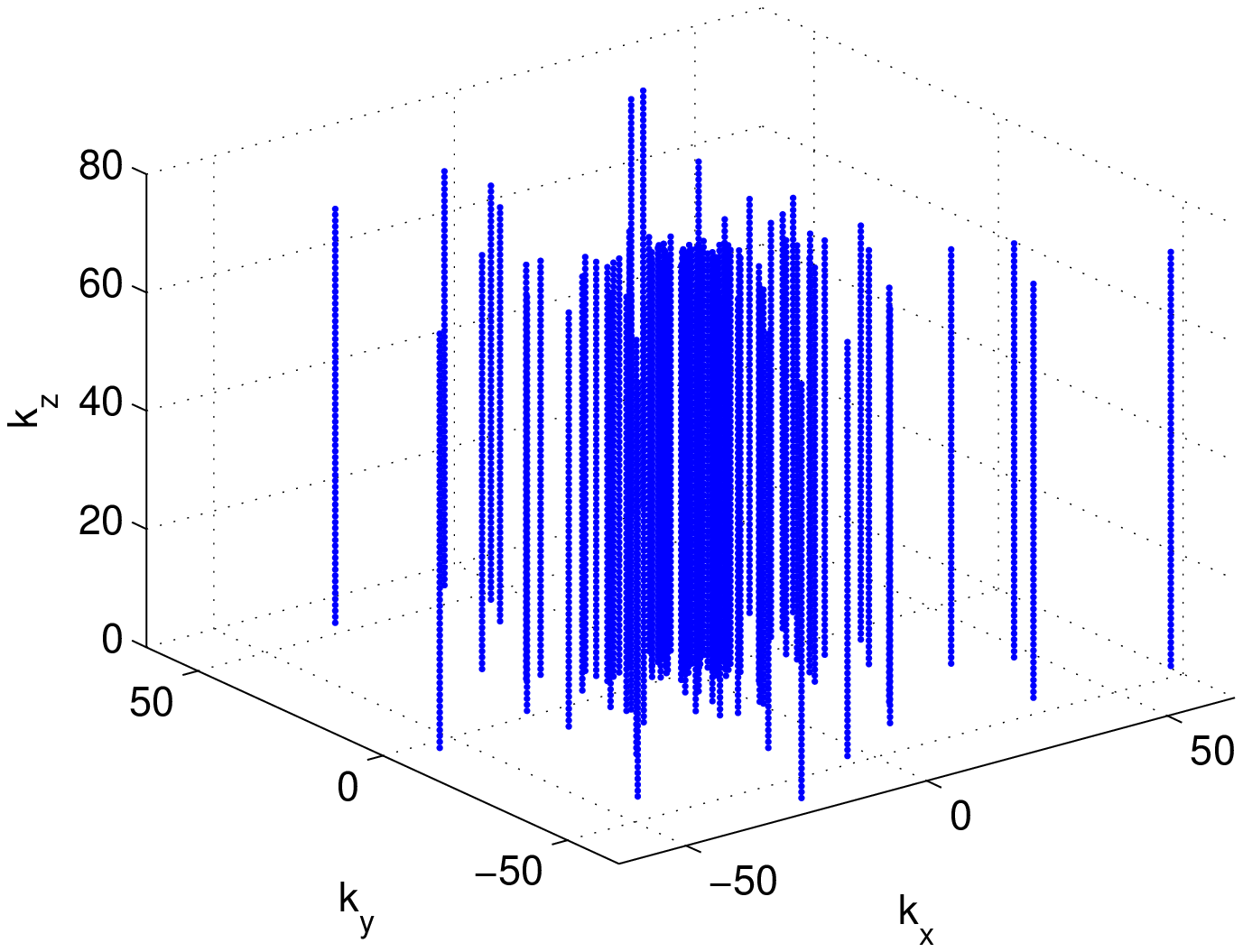}\label{kpoints3Dpois}}&
\subfigure[]{\includegraphics[width=3in]{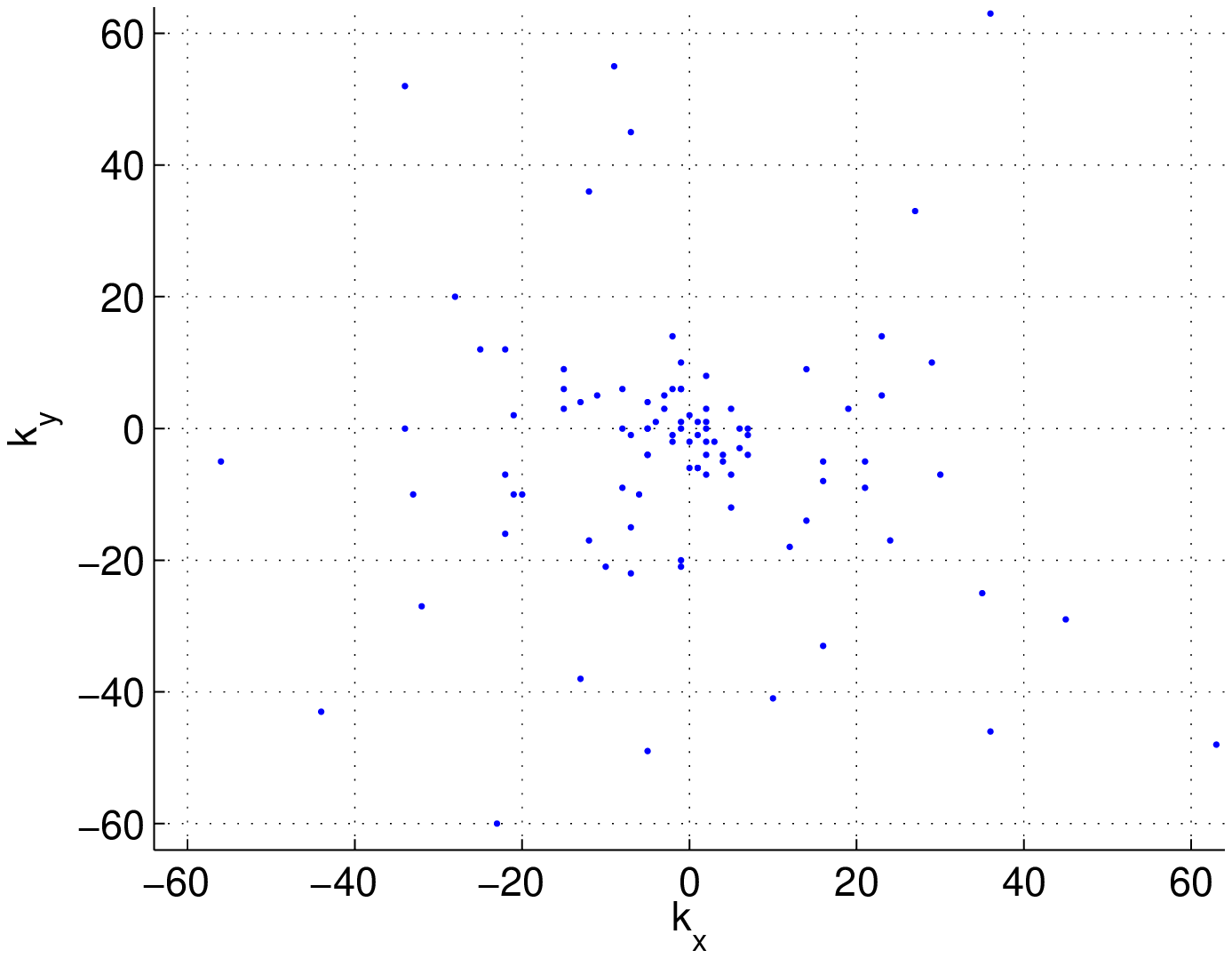}\label{kpointsxypois}}
\end{tabular}
\vspace{-0.1in}
\caption{Volumetric Image Reconstruction: Acquired points for data observation (based on variable-density Poisson disk undersampling) a) 3-dimensional view b) projection of the points on $k_x-k_y$ plane. Readout lines are the lines passing through these 100 points.}\label{3Dkpoints_poisson}
\vspace{-0.1in}
\end{figure*}

To validate performance of the proposed undersampling technique, we also implemented the approach with two other data acquisition schemes. The first one is rectilinear undersampling, which can be obtained using multi-shot Echo Planar Imaging (EPI) and is widely used in clinical applications. We acquired 100 $k-$space lines, uniformly distributed over $k_x-k_y$ plane, as shown in \Fig{rect3d}. Variable-density Poisson disk undersampling \cite{mitchell2012variable} was also used as another scheme for data acquisition. These points are illustrated in \Fig{3Dkpoints_poisson}. In order to be consistent in our comparisons, we used variable-density Poisson disk sampling and rectilinear approach to generate 100 points over $k_x - k_y$ plane and then all the points on the lines passing through these points were considered for data acquisition. 

\subsection*{Error Analysis}
\begin{figure*}[htb!]
\centering
\begin{tabular}{ccc}
\hspace{-0.2in}\subfigure[]{\includegraphics[width=2.5in]{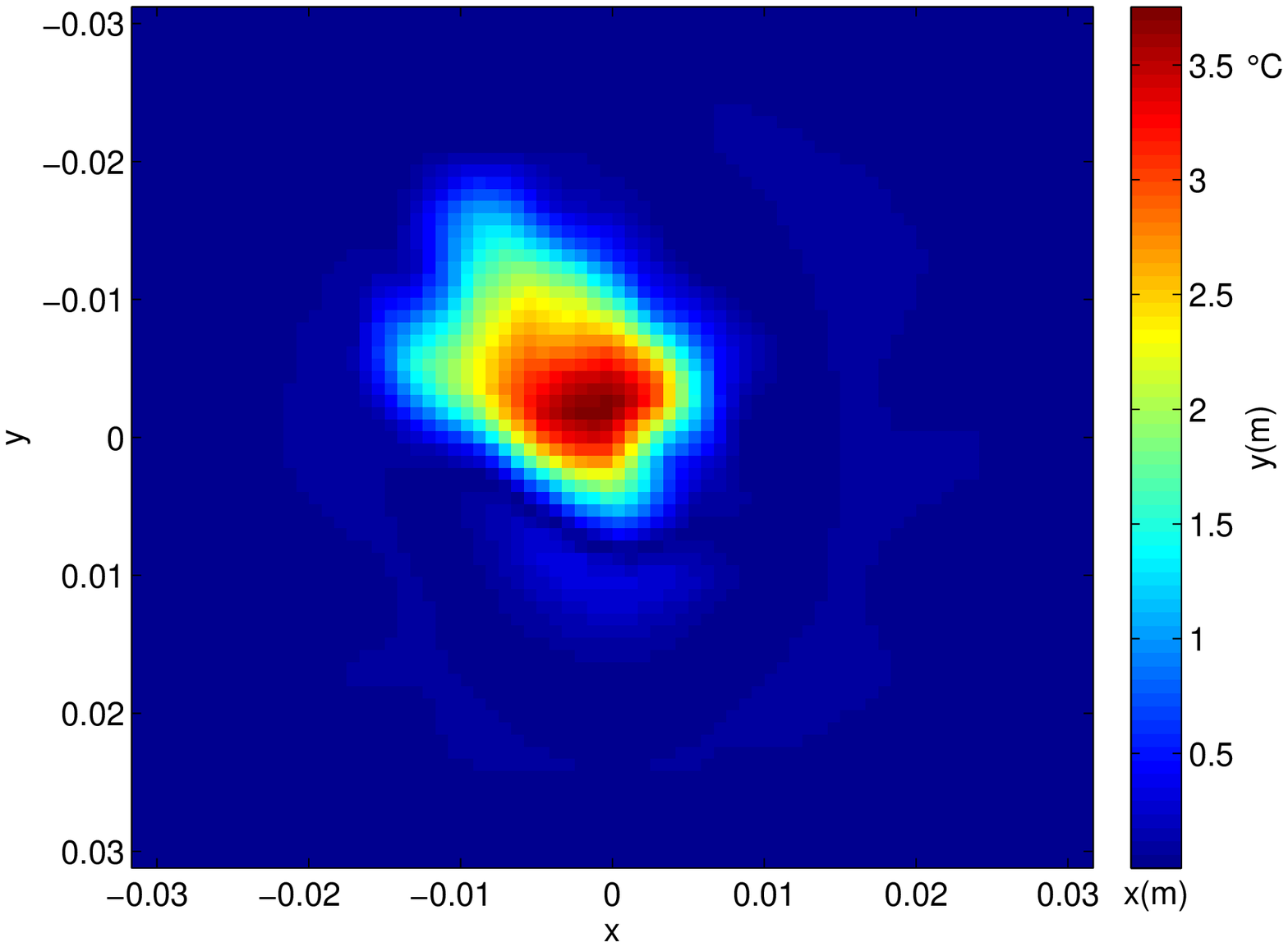}\label{err40_var}} &
\hspace{-0.2in}\subfigure[]{{\includegraphics[width=2.5in]{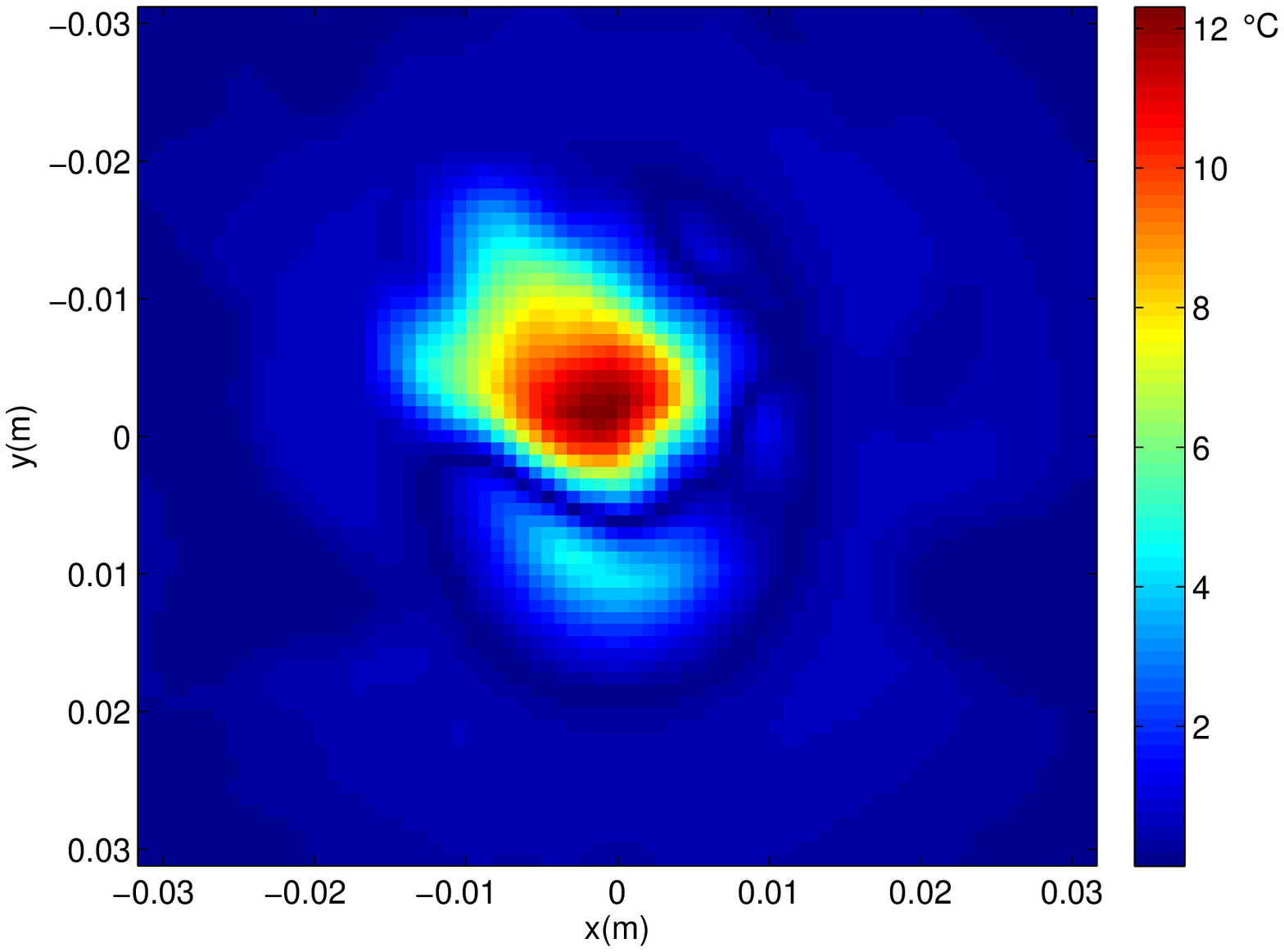}}\label{err40_rect}} &
\hspace{-0.2in}\subfigure[]{\includegraphics[width=2.5in]{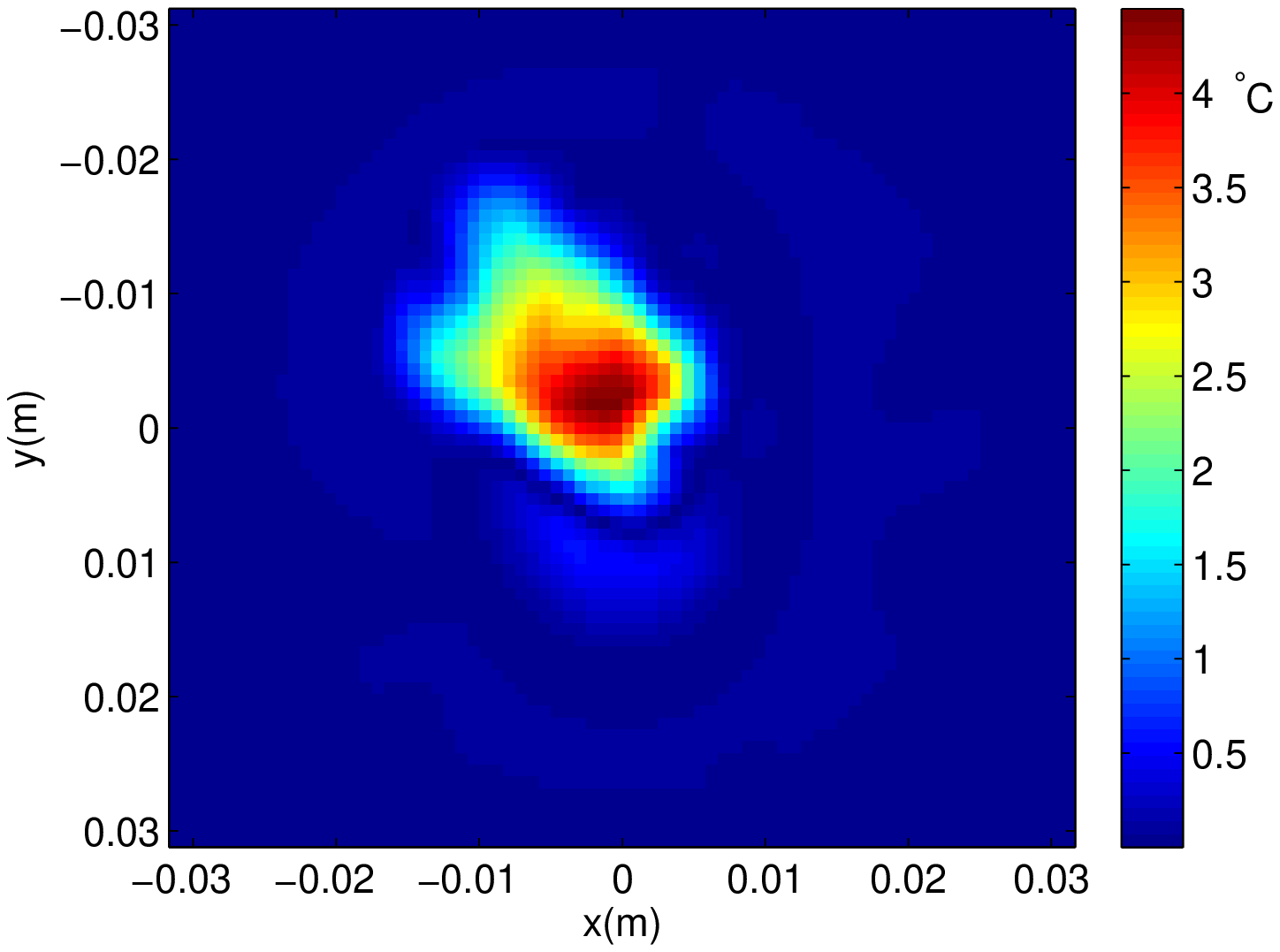}\label{err40_pois}}\\
\end{tabular}
\caption{Volumetric Image Reconstruction: The error between posterior
estimate and actual temperature field over the $40^{th}$ slice, obtained based on a) Proposed technique, b) Rectilinear undersampling, c) Variable-density Poisson disk undersampling}\label{err40} 
\end{figure*}

Obtained posterior estimates of $\mu$ from each of these data acquisition techniques are then used in the Pennes bioheat equation to estimate the
true temperature field. {\color{black} \Fig{err40} represents the error between the true temperature field and its estimate in a representative  $x-y$ plane, obtained by different data acquisition schemes.}
As illustrated in \Fig{err40},
when the variance maximization is used for data acquisition, the maximum discrepancy
between the estimated and true temperature field is less than $4\ ^\circ$C, while both other undersampling schemes result in higher value of maximum error.

\begin{figure}[htb!!!]
\centering
\vspace{-0.1in}
{\includegraphics[width=3.35in,height=2.4in]{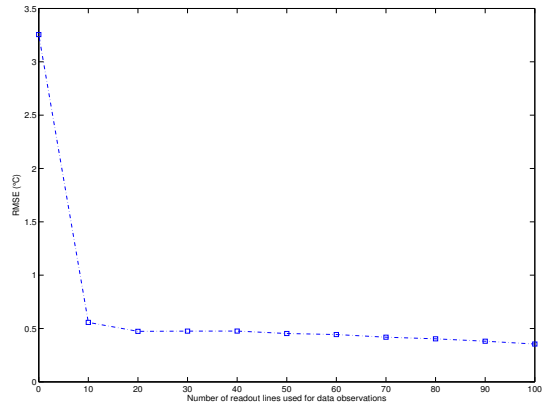}}
\vspace{-0.1in}
\caption{Volumetric Image Reconstruction: RMSE (over all voxels in the full 3D volume) between the estimated and actual temperature fields versus different number of observations (acquired based on the proposed approach).}\label{errcvr} 
\end{figure}

To further study performance of the proposed technique, we have shown the convergence of the Root Mean Square Error (RMSE) between the temperature estimate and its actual value while using the proposed technique for data acquisition. As shown in \Fig{errcvr}, the RMSE in the 
temperature field estimate consistently reduces with increasing the number of readout lines. 

To summerize, Table \ref{table_err3D} represents the maximum and the RMSE between the true temperature field and its estimate, calculated over all voxels in the full 3D volume, when maximum heating occurs. Note that without any observation data the maximum difference between the actual temperature map and its estimate is more than $33\ ^\circ$C and by increasing the number of readout lines, it consistently reduces until it reaches to almost $4\ ^\circ$C when 100 readout lines (obtained using the proposed data acquisition technique) are used. Also, the RMSE between the temperature estimate and its
actual value reduces by a factor of 9.  Table \ref{table_err3D} also shows the maximum error and RMSE while using rectilinear undersampling and variable-density Poisson disk undersampling. It is clear that variable-density Poisson disk undersampling performs better than rectilinear undersampling, but comparing these values with corresponding error metrics while using variance maximization technique clearly shows performance of the proposed subsampling technique.
\begin{table*}[htb!!!]
\centering
\footnotesize
\caption{Volumetric Image Reconstruction: Convergence of error (calculated over all voxels in the full 3D volume) given
different numbers of $k-$space samples for model-data fusion. Prior values of error are shown in second column (with zero observations).}\label{table_err3D}
\begin{tabular}{lccccc}
\hline
Number of readout lines & 0 & 30 & 50 & 80 & 100 \\
\hline
\hline
maximum error (Proposed Technique) $^\circ$C & 33.08 & 5.16 & 4.98 & 4.37 & 3.81  \\
\hline
RMSE (Proposed Technique) $^\circ$C  & 3.26 & 0.47 & 0.45 & 0.40 & 0.35 \\
\hline
\hline
\multicolumn{1}{c}{maximum error (Rectilinear Undersampling) $^\circ$C}  & 33.08 & & & & 12.53  \\
\hline
 \multicolumn{1}{c}{RMSE (Rectilinear Undersampling) $^\circ$C } &  3.26 & & & & 1.23 \\
\hline
\hline
\multicolumn{1}{c}{maximum error (Poisson Undersampling) $^\circ$C}  & 33.08 & & & & 4.52  \\
\hline
 \multicolumn{1}{c}{RMSE (Poisson Undersampling) $^\circ$C } &  3.26 & & & & 0.41 \\
\hline
\hline
\end{tabular}
\vspace{0.05in}
\vspace{-0.2in}
\end{table*}


\subsection{Planar Image Reconstruction}
\begin{figure}[htb!!]
\vspace{-0.1in}
\centering
\includegraphics[width=3.5in,height=2.35in]{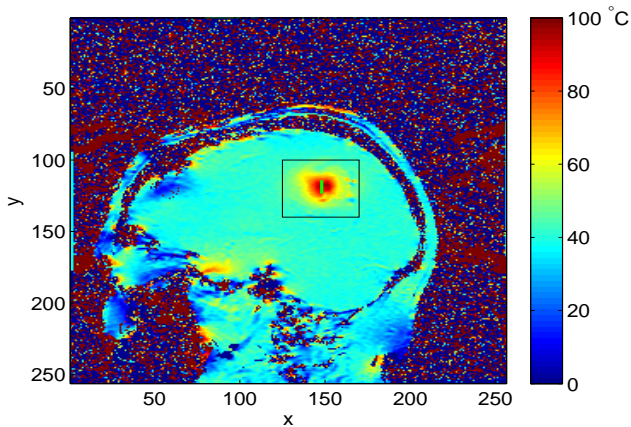}
\vspace{-0.1in}
\caption{Planar Image Reconstruction: Position of the tumor and the
temperature field resulting from MRgLITT (FOV = 26 cm$^2$). Black box and green line inside represent the region of interest and laser fiber, respectively.}\label{layout2D}
\vspace{-0.05in}
\end{figure}

{\color{black}
Performance of the proposed subsampling approach is retrospectively
validated in fully sampled
planar temperature imaging acquired during the thermal ablation process
\textit{in vivo} human brain. 
}
The LITT procedure (power = 11.85 watts for a period of 94 sec.) was monitored 
in real-time using the temperature-sensitive 
PRF shift technique acquired with a 2D spoiled gradient-echo to
generate temperature measurements using GE scanner at every $\Delta t
=$5 seconds (TR/FA = 38 ms/30$^\circ$, frequency $\times$ phase
= 256 $\times$ 128, FOV = 26 cm$^2$, BW = 100 Hz/pixel,
slice thickness 5 mm).
 \Fig{layout2D} shows position of the tumor and
the temperature resulting from laser irradiation in the region.
\eq{obsmodel} is used to model the MR signal.  Numerical values of the
parameters involved in sensor structure are shown in Table
\ref{sensor_para_human2D}. As shown in Table \ref{sensor_para_human2D}, we assumed $T_1$ to be a linear function of tissue temperature \cite{rieke2008mr} in order to incorporate the possible changes in its value due to heating. The value of $T_1$ at body temperature, denoted by $T_1^0$ is considered to be $1.05$ sec.. Also, a normalized MRI image (without any
heating involved) is used to incorporate the effect of inhomogeneities in
magnetization $M(\x)$. 

\begin{figure}
\centering
\vspace{-0.1in}
\begin{tabular}{c}
\subfigure[]{\includegraphics[width=3.2in,height=2.15in]{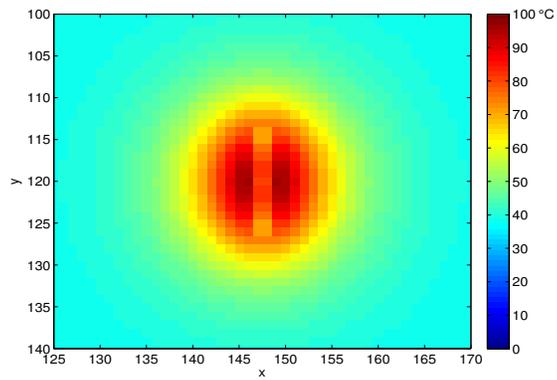}\label{mean_cvr_human_2d}} \\
\hspace{0in}\subfigure[]{\includegraphics[width=3.2in,height=2.15in]{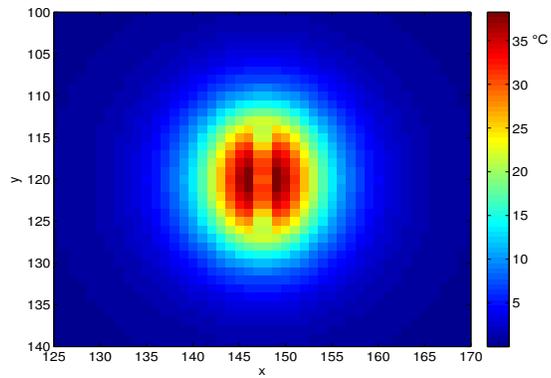}} 
\end{tabular}
\caption{
Planar Image Reconstruction: statistics of temperature field over ROI (ROI = 55 mm $\times$ 60 mm), resulted from uncertain $\mu$, a) mean ($^\circ$C), b) standard deviation ($^\circ$C).}\label{Tstat2d}
\vspace{-0.15in}
\end{figure}

{ \color{black}
Here, we assume that the optical attenuation
coefficient $\mu$ is spatially homogeneous and uniformly distributed between
$100\ \frac{1}{\text{m}}$ and $400\ \frac{1}{\text{m}}$, \textit{i.e.} $\mu\sim \mathcal{U}[100,400]$. 
Note that
$\mu$ is the only uncertain parameter while solving the Pennes bioheat
equation; all the other parameters are assumed to be known from the
physical properties of the tissue,
Table \ref{pennes_para}.  }
A set of 15 Gauss-Legendre quadrature points \cite{hildebrand1987introduction} (which accurately integrate up to the $29^{th}$ order polynomials), spread over the range of uncertain parameter $\mu$, are used to
quantify the uncertainty in the temperature field. \Fig{Tstat2d} illustrates
mean and standard deviation of the associated temperature field within the region
of interest (ROI).

\begin{table}[htb!]
\vspace{-0.1in}
\centering
\caption{Planar Image Reconstruction: value of the parameters involved in MRI signal model}\label{sensor_para_human2D}
\begin{tabular}{cc}
Parameter & Value\\
\hline
$\theta$  & $\pi/6$ rad \\
\hline
$T2*$ & $70$ ms \\
\hline
$T_E$ & $20$ ms\\
\hline
$\gamma$ & 42.58 MHz/T\\
\hline
$\alpha$ & $-0.0102$ ppm/C\\
\hline
$B_0$ & $1.5$ T\\
\hline
$T_1$ & $T_1=T_1^0(1+0.01\Delta u),\quad T_1^0=1.05$ ms\\
\hline
\end{tabular}
\vspace{-0.in}
\end{table}

The proposed technique in Section \ref{sec:dyn_data} is utilized to approximate the
observation points on $k-$space with largest information content.
Given uncertainties in the model, \eq{pennes},
\Fig{human2dpoints} shows 20 readout lines that are predicted to have the highest
variance in data acquisition.
Note that as expected, readout lines pass through the points in center of the
$k-$space, i.e. the points with the greatest information content. This corroborates the fact that the points in center
of the $k-$space contain more information regarding the main features of the
image. 

The predicted $k-$space measurements are then retrospectively 
extracted from the fully sampled data and used in a minimum variance
framework to estimate posterior statistics of parameter $\mu$, as described
in Section \ref{sec:gpc_minvariance}. We considered SNR = 25 while using minimum variance framework. \Fig{human2d_cvr} demonstrates
convergence behavior for the mean and variance of $\mu$ versus different
number of observations. \textcolor{black}{Note that only $\frac{20}{256}\simeq 0.08$ of
the total number of readout lines in $k-$space are used for estimation of 
$\mu$. This results in a speed up factor of 12 (based on division of number of phase lines to the number of readout lines) in comparison with full-sampling of $k-$space, \Fig{human2dpoints}}.
As a reference, a global optimization
to find the value of $\mu$ which results in the least discrepancy between
the predicted and fully sampled temperature measurement data was performed. 
This is shown with a dashed black line in \Fig{human2d_cvr}. 
Results obtained by the proposed approach converges with the
result obtained from global optimal and fully sampled technique.

\begin{figure}
\centering
\vspace{-0.1in}
\includegraphics[width=3.2in]{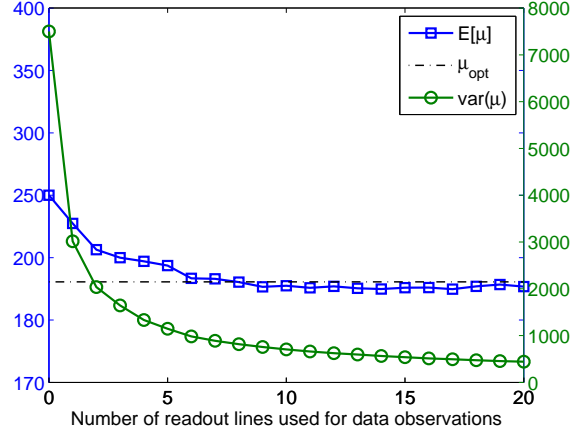}
\caption{Planar Image Reconstruction: Convergence of $\mu$ estimate versus number of readout lines. Square-blue and circle-green lines represent mean and variance of $\mu$ estimate, respectively. Dashed black line represents the optimal value of $\mu$, obtained from global optimization.}\label{human2d_cvr}\label{human2d_cvr}
\vspace{-0.1in}
\end{figure}

\begin{figure}[htb!!]
\centering
\vspace{-0in}
{\includegraphics[width=3.4in]{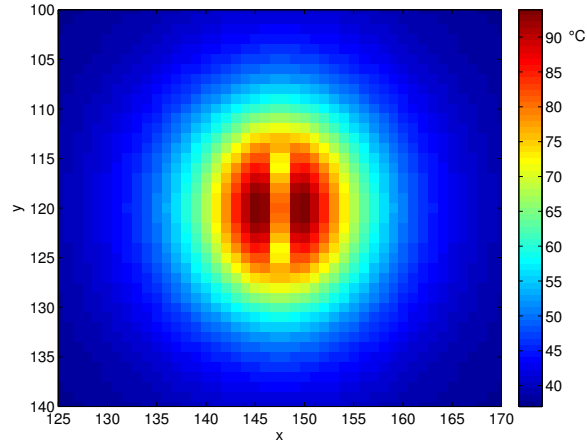}}
\vspace{-0.1in}
\caption{Planar Image Reconstruction: Posterior estimate of the temperature map (ROI = 55 mm $\times$ 60 mm). }\label{human_Tpost}
\vspace{-0.1in}
\end{figure}

The posterior estimate of $\mu$ is then used in Pennes bioheat equation to estimate the temperature field of the tissue. \Fig{human_Tpost} illustrates the posterior estimate of the temperature field.

\subsection*{Comparison with Other Undersampling Schemes}
\begin{figure*}[htb!!]
\centering
\begin{tabular}{ccc}
\hspace{0in}\subfigure[]{\includegraphics[width=2.1in]{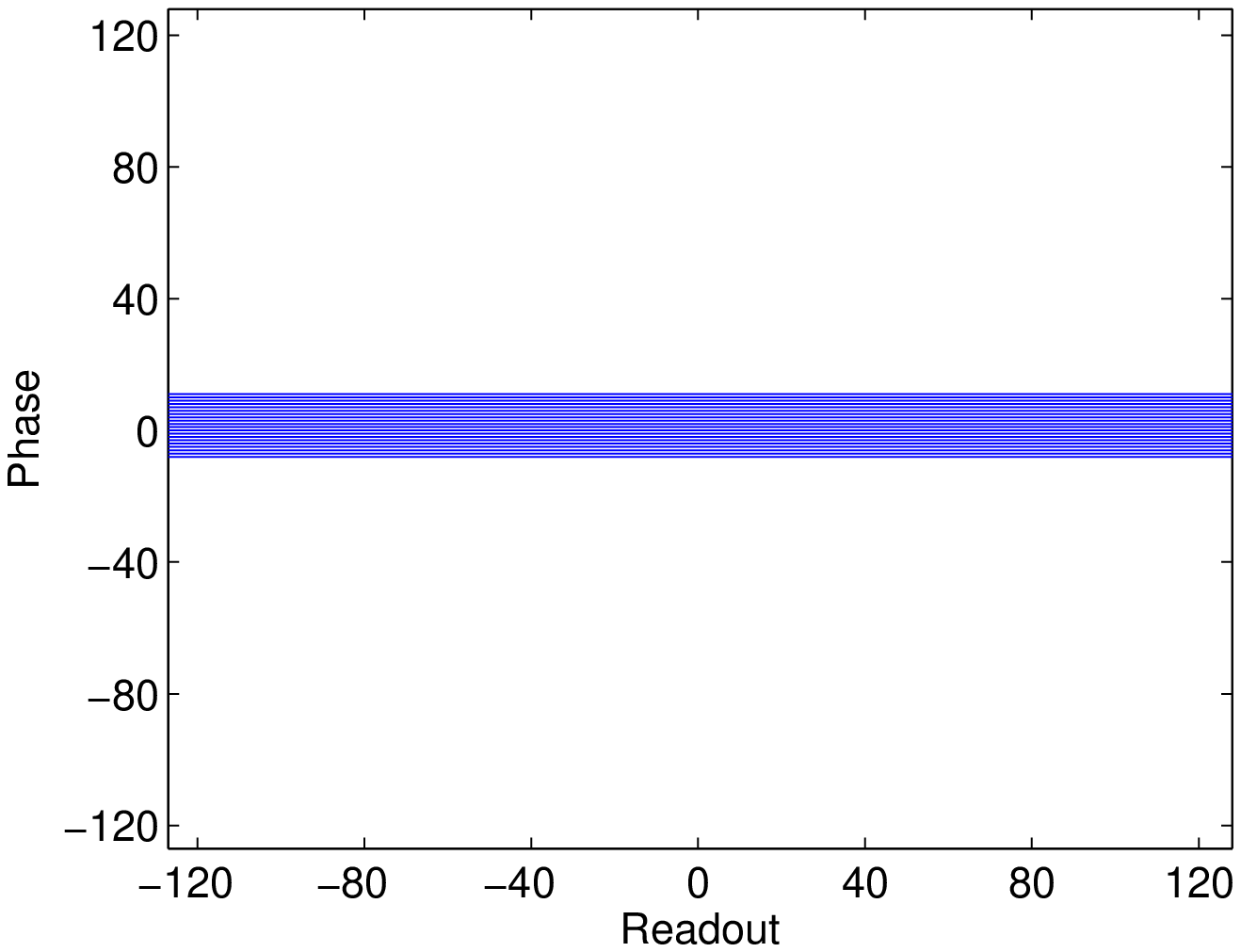}\label{human2dpoints}} &
\hspace{0in}\subfigure[]{\includegraphics[width=2.1in]{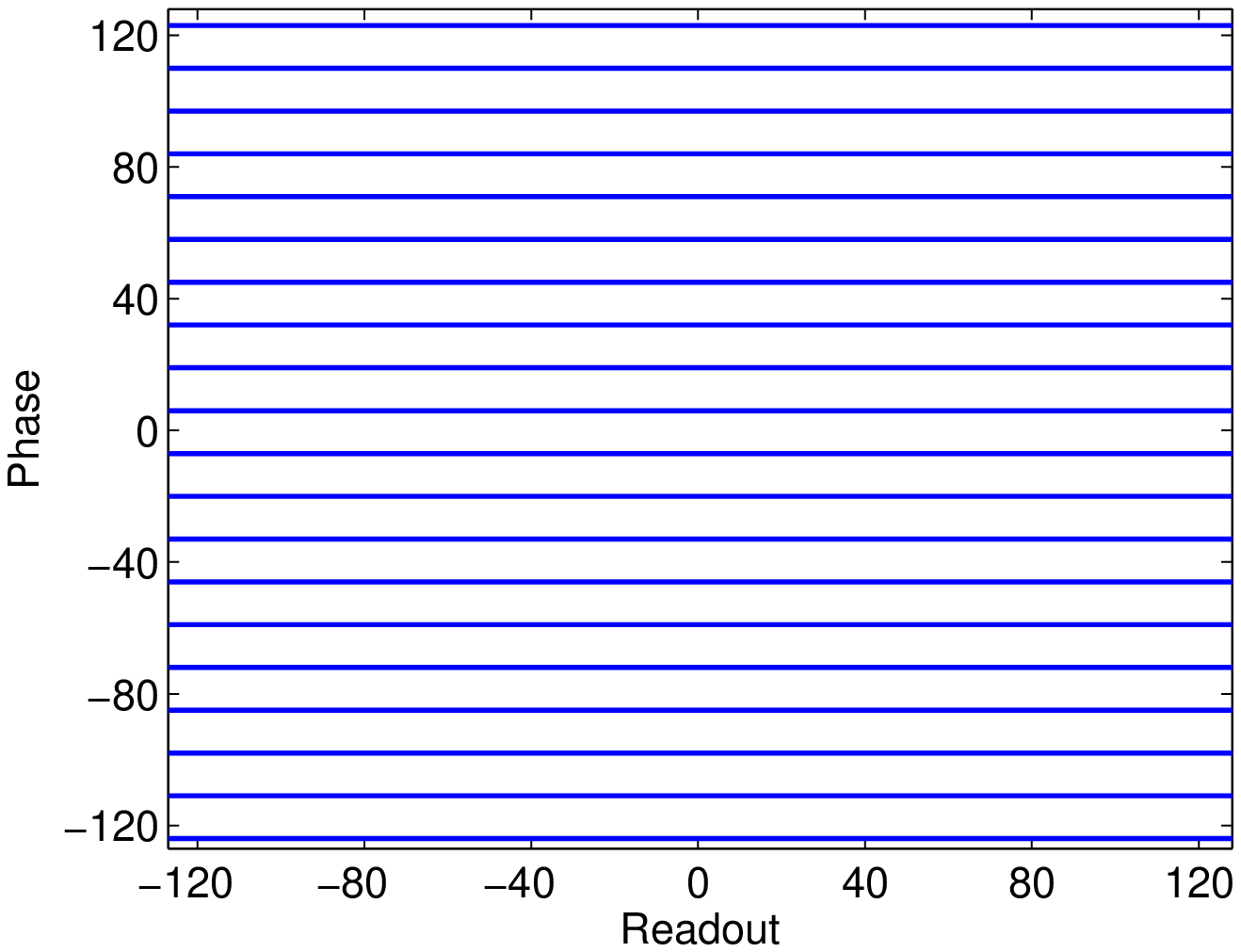}\label{human2dpoints_uniform}} &
\hspace{0in}\subfigure[]{\includegraphics[width=2.1in]{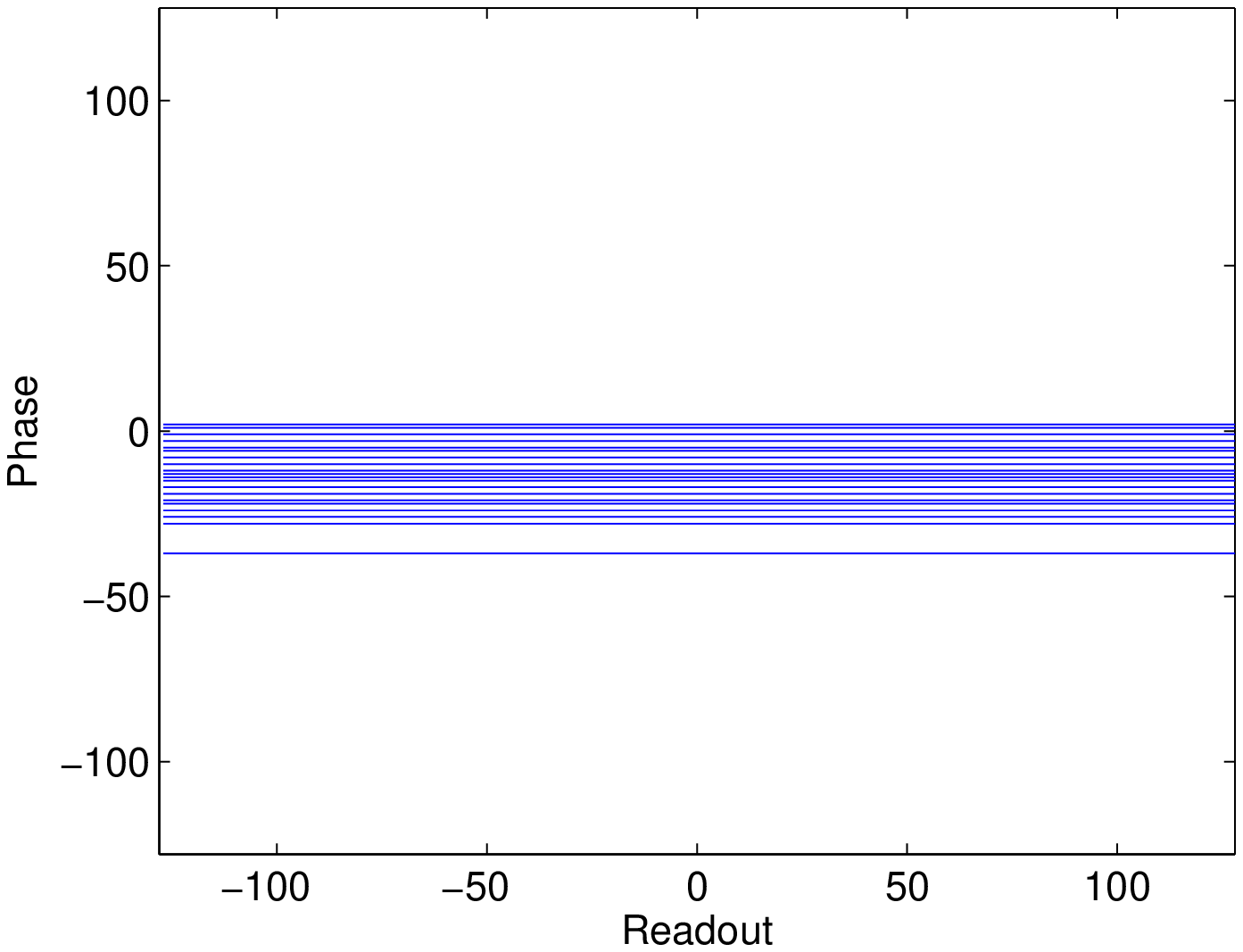}\label{human2dpoints_poisson}}
\end{tabular}
\vspace{-0.1in}
\caption{Planar Image Reconstruction: a) The 20 readout lines used for data acquisition, obtained from variance maximization. b) Uniformly distributed readout lines based on rectilinear undersampling. c) Readout lines based on variable-density Poisson disk undersampling. }
\vspace{-0.1in}
\end{figure*}

Similar to previous example, we have also used rectilinear and variable-density Poisson undersampling techniques for $k-$space data acquisition. \Fig{human2dpoints_uniform} and \Fig{human2dpoints_poisson} show 20 readout lines, generated based on rectilinear and variable-density Poisson disk undersampling techniques, respectively.

{\color{black}
Given each data acquisition scheme, the predicted $k-$space measurements are used in minimum variance
framework to estimate posterior statistics of parameter $\mu$. The obtained posterior estimate of $\mu$ is then used in Pennes
bioheat equation to estimate the temperature field of the tissue.




\subsection*{Error Analysis}
The error (calculated over ROI) between the posterior estimate of the temperature field and actual
temperature map is illustrated in \Fig{human2d_err}. It is clear from \Fig{human2d_err} that the proposed technique results in less error comparing with both rectilinear and variable-density Poisson disk undersampling. Note that when using the proposed method most of error occurs in the edge of ROI due to presence of artifacts, while a considerable amount of error occurs in the center of ROI while using rectilinear or Poisson disk undersampling . 
\begin{figure*}[htb!!]
\centering
\vspace{-0in}
\begin{tabular}{ccc}
\hspace{-0.3in}{\includegraphics[width=2.6in]{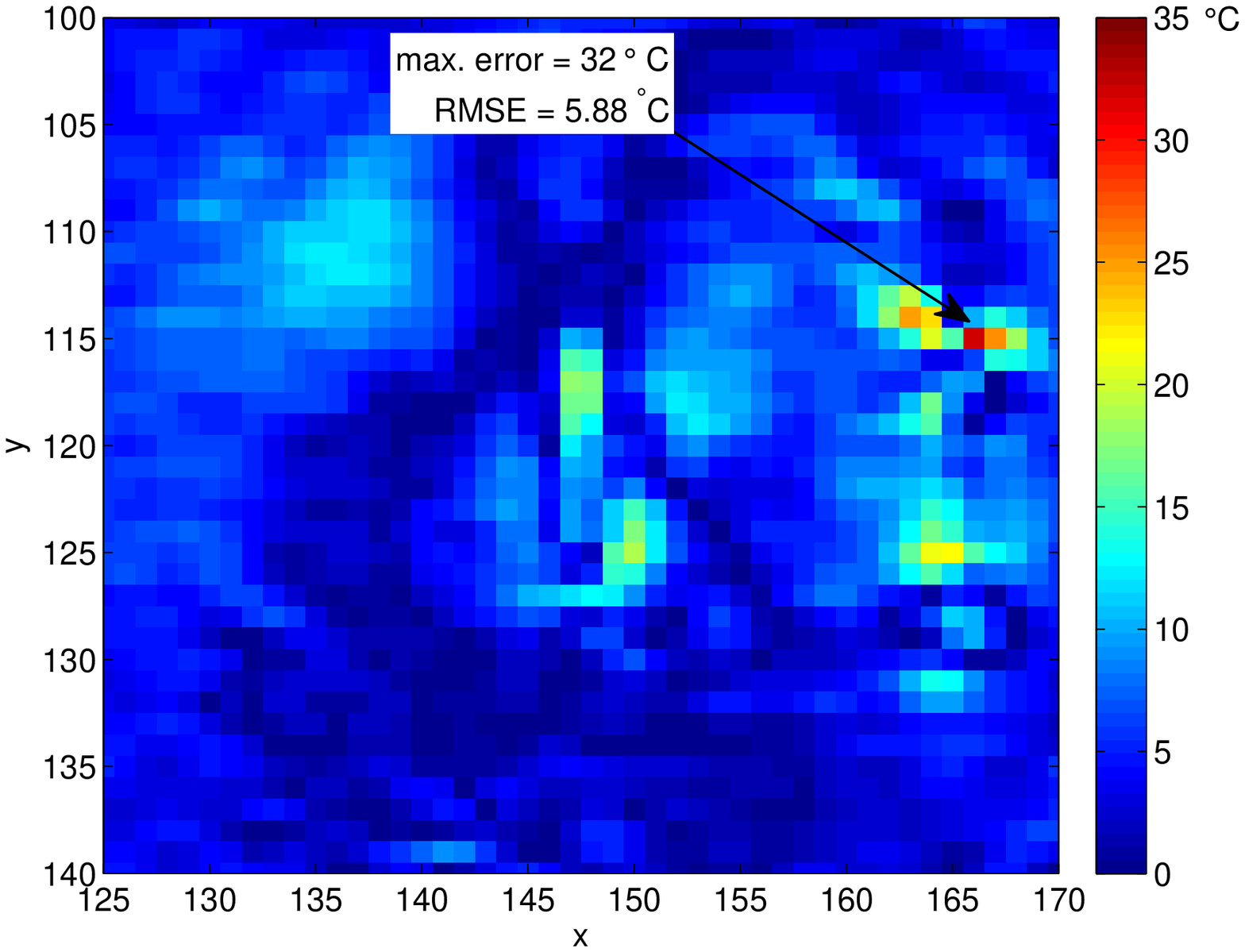}} &
\hspace{-0.3in}{\includegraphics[width=2.6in]{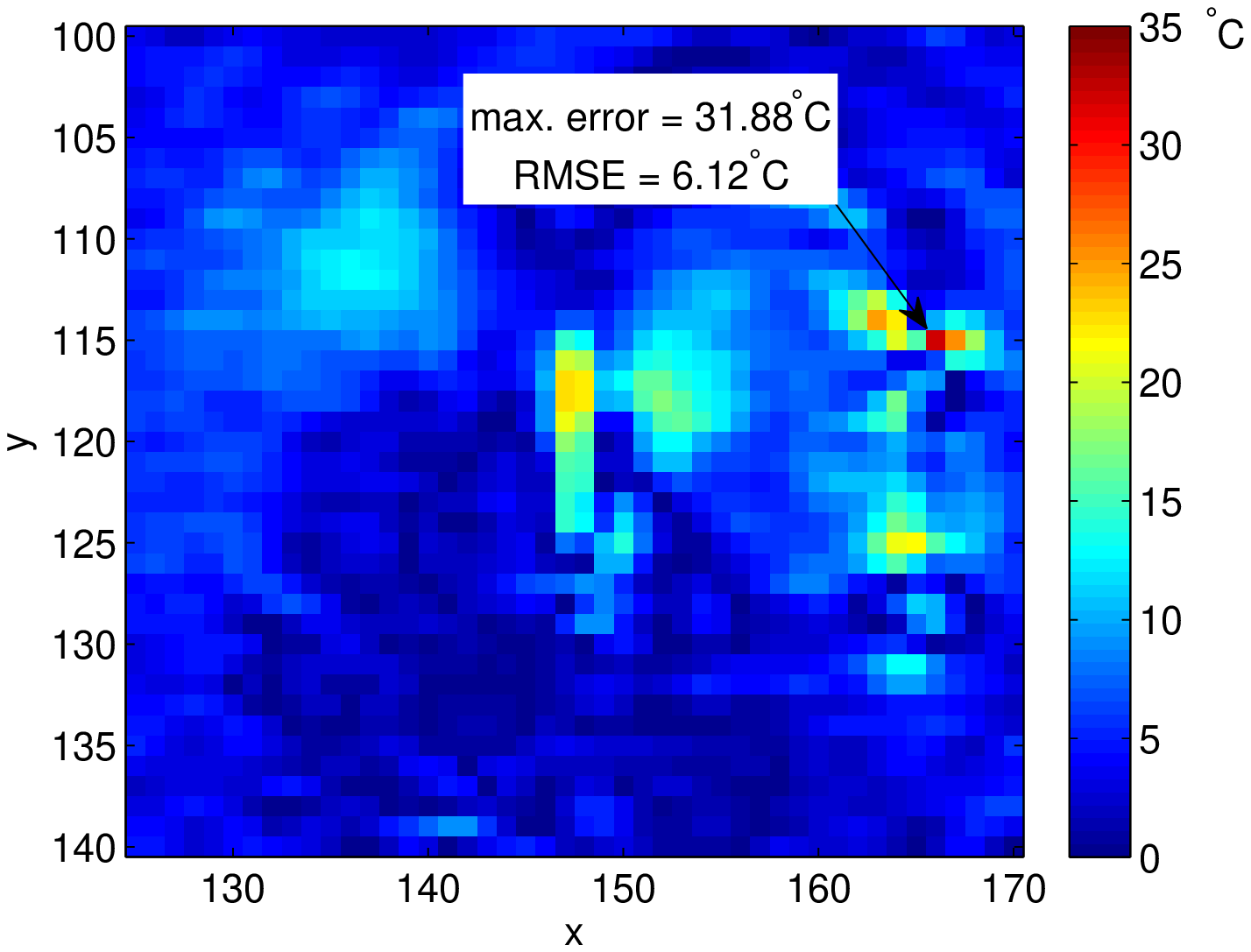}} &
\hspace{-0.3in}{\includegraphics[width=2.6in]{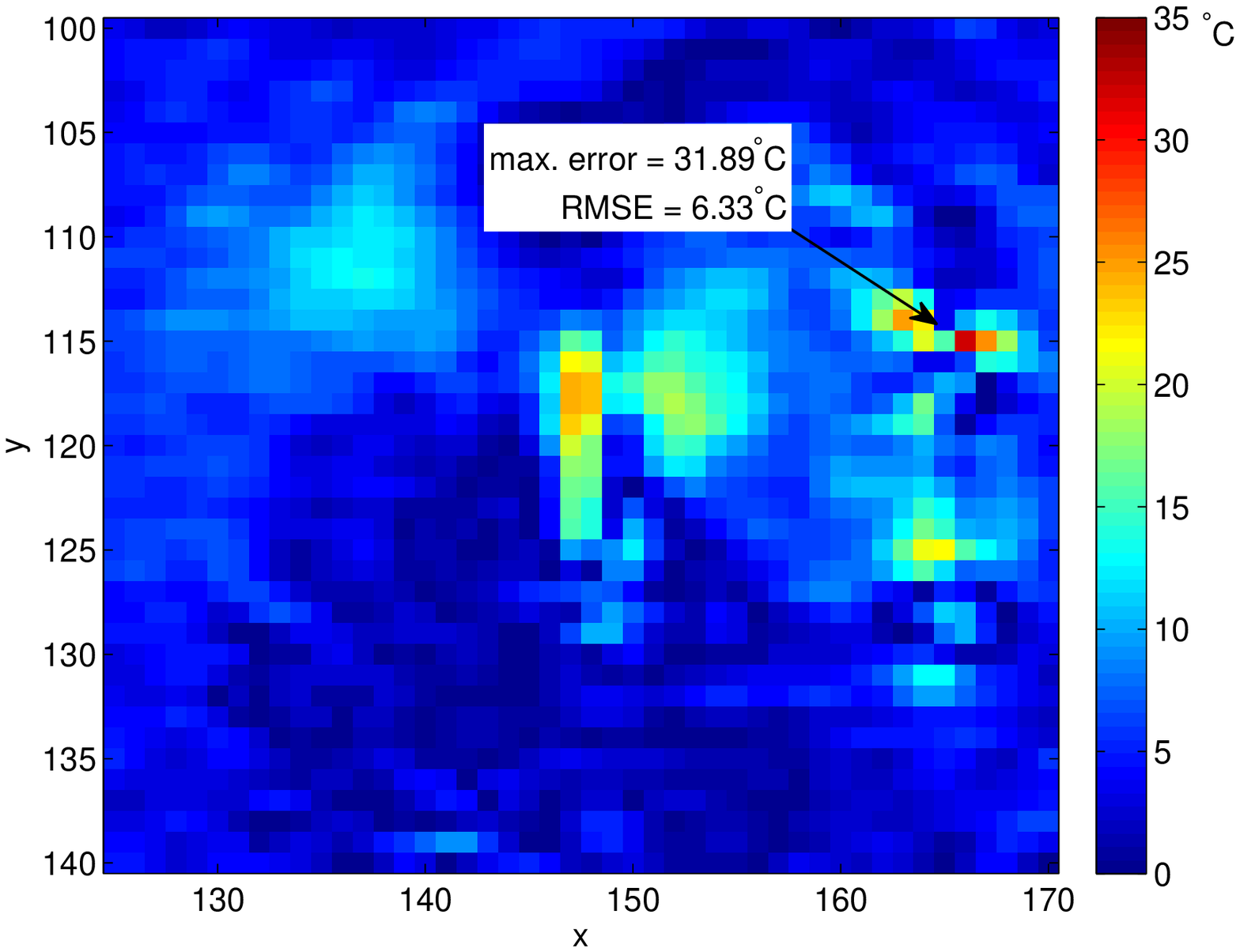}} 
\end{tabular}
\vspace{-0.1in}
\caption{Planar Image Reconstruction: a) The error between the posterior estimate  and actual measurement of temperature field, obtained by the proposed technique. Note that the maximum error occurs over the regions close to border. The high value of error in regions near the border of the figure is due to the presence of artifacts in MR temperature data. b) The error between the posterior estimate  and actual measurement of temperature field, obtained by using rectilinear undersampling. c) The error between the posterior estimate  and actual measurement of temperature field, obtained by using variable-disk Poisson disk undersampling. Note the significant amount of error in center of ROI for cases a and b (ROI = 55 mm $\times$ 60 mm). }\label{human2d_err}
\vspace{-0.1in}
\end{figure*}

To study the effect of the number of data observations on performance of the proposed technique, we have shown the RMSE between the estimated and fully-sampled temperature field in \Fig{human2derrcvr}, where in general the value of RMSE reduces with increasing the number of readout lines. Note that a major portion of RMSE is due to presence of artifacts.
\begin{figure}[htb!!]
\centering
\includegraphics[width=3.25in]{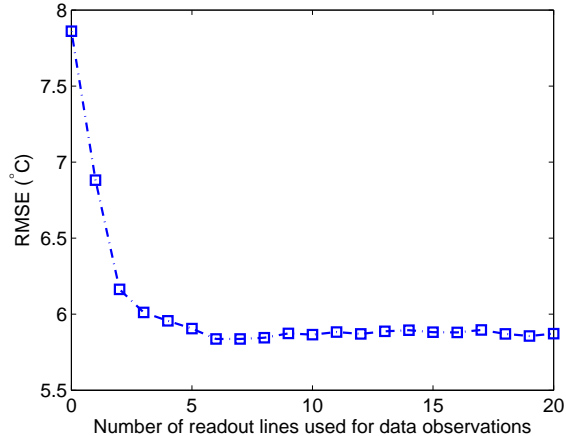}
\caption{Planar Image Reconstruction: Convergence of RMSE (calculated over ROI) between the acquired temperature map from fully-sampled $k-$space and temperature estimate (obtained from the proposed technique). Note that a major portion of RMSE is due to the presence of artifacts.}\label{human2derrcvr}
\end{figure}

\subsection{Phantom Experiment}
To further validate the proposed approach, we designed a phantom experiment for 3 dimensional temperature imaging. A tissue-mimicking agar phantom ($1.5\%$ gel by weight) was used to simulate the thermal ablation process. The LITT procedure (power = 1 watt for a period of 10 min.) was monitored 
in real-time using the temperature-sensitive 
PRF shift technique acquired with EFGRE3D to
generate temperature measurements using a GE scanner at every $\Delta t
=$15 seconds (TR/FA = 5 ms/5$^\circ$, frequency $\times$ phase
= 256 $\times$ 256 $\times$ 26, FOV = 12.8 cm $\times$ 12.8 cm $\times$ 26 mm , BW = 390.62 Hz/pixel,
slice thickness 1 mm). Temperature images were recorded over a number of 26 slices in $z-$direction (each with thickness of 1 mm). \Fig{phantom} shows the position of the laser fiber and corresponding ROI. The parameters used in phantom experiment are shown in Tables \ref{phantom_para1} and \ref{phantom_para2}. Note that we assumed a linear dependence between $T_1$ and temperature in order to incorporate possible changes of $T_1$ due to changes in temperature  \cite{vesanen2013temperature}. Also, a normalized MRI image (without any heating involved) was used to account for the effect of inhomogeneities in magnetization $M(\x)$. 
\begin{figure}[htb!!]
\centering
\includegraphics[width=3.5in]{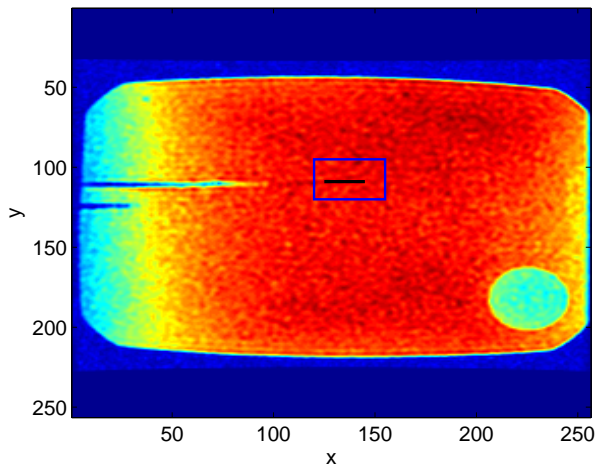}
\caption{Phantom Experiment: blue rectangle and black line inside represent the ROI and position of laser fiber, respectively.}\label{phantom}
\end{figure}

\begin{table}[htb!!]
\caption{
\color{black}
\textcolor{black}{Phantom experiment: Data used in numerical simulations ~\cite{taylor2008dynamic}}}\label{phantom_para1}
\centering
\begin{tabular}{ccccc} \hline
$\Lambda $ $ \frac{W}{ m \cdot K}$ & $\omega$ $\frac{kg}{m^3 s}$ &  $\rho$
$\frac{kg}{m^3}$ &   $c_{blood}$ $ \frac{J}{kg \cdot K}$ &  $c$
$\frac{J}{kg \cdot K}$ \\ \hline\hline
          0.6              &             0           &  1000
&           0                     &                  3900          \\
\hline
\end{tabular}
\end{table}

\begin{table}[htb!]
\centering
\caption{Phantom experiment: parameter values involved in sensor model}\label{phantom_para2}
\begin{tabular}{cc}
Parameter & Value\\
\hline
\hline
$\theta$ (degree) & $5$ \\
\hline
$T_1$ (sec.) & $2.56+0.1\Delta u$ \\
\hline
$T_2^*$ (ms) & 30\\
\hline
$\Delta \omega_0$ (rad) & 0\\
\hline
$T_E$ (ms) & $1.676$ \\
\hline
$T_R$ (ms) & $5$\\
\hline
$\gamma$ (MHz/T) & 42.58 \\
\hline
$\alpha$ (ppm/C)& $-0.0102$ \\
\hline
$B_0$ (T)& $3$ \\
\hline
$u_0$ ($^\circ$C)& $19$\\
\hline
\end{tabular}
\end{table}

We considered the only uncertain parameter to be the optical attenuation coefficient, which is assumed to be uniformly distributed between 0 and $200$ $\frac{1}{\text{m}}$, i.e. $\mu \in \mathcal{U}(0,200)$. We used a set of 11 Gauss-Legendre quadrature points \cite{hildebrand1987introduction} (which accurately integrates up to the $21^{st}$ order polynomials) to quantify the effect of uncertain parameter $\mu$ on the temperature field and $k-$space signal. Then, the proposed technique was utilized to approximate $k-$space locations with highest amount of variance. \Fig{phantomkpoints} illustrates the $k-$space points found by this approach.

\begin{figure}[htb!!!]
\centering
\begin{tabular}{c}
\subfigure[]{\includegraphics[width=3.2in,height=2.2in]{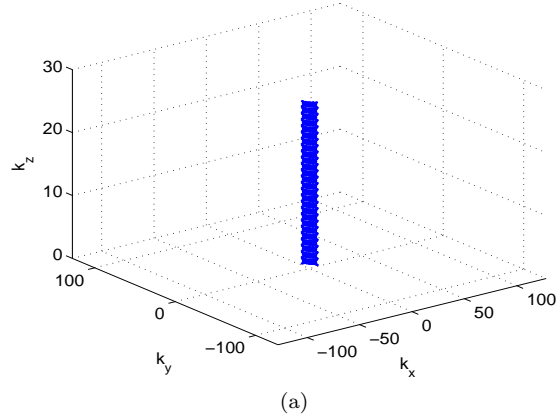}\label{phantom3dpoints}} \\
\subfigure[]{\includegraphics[width=3.2in,height=2.2in]{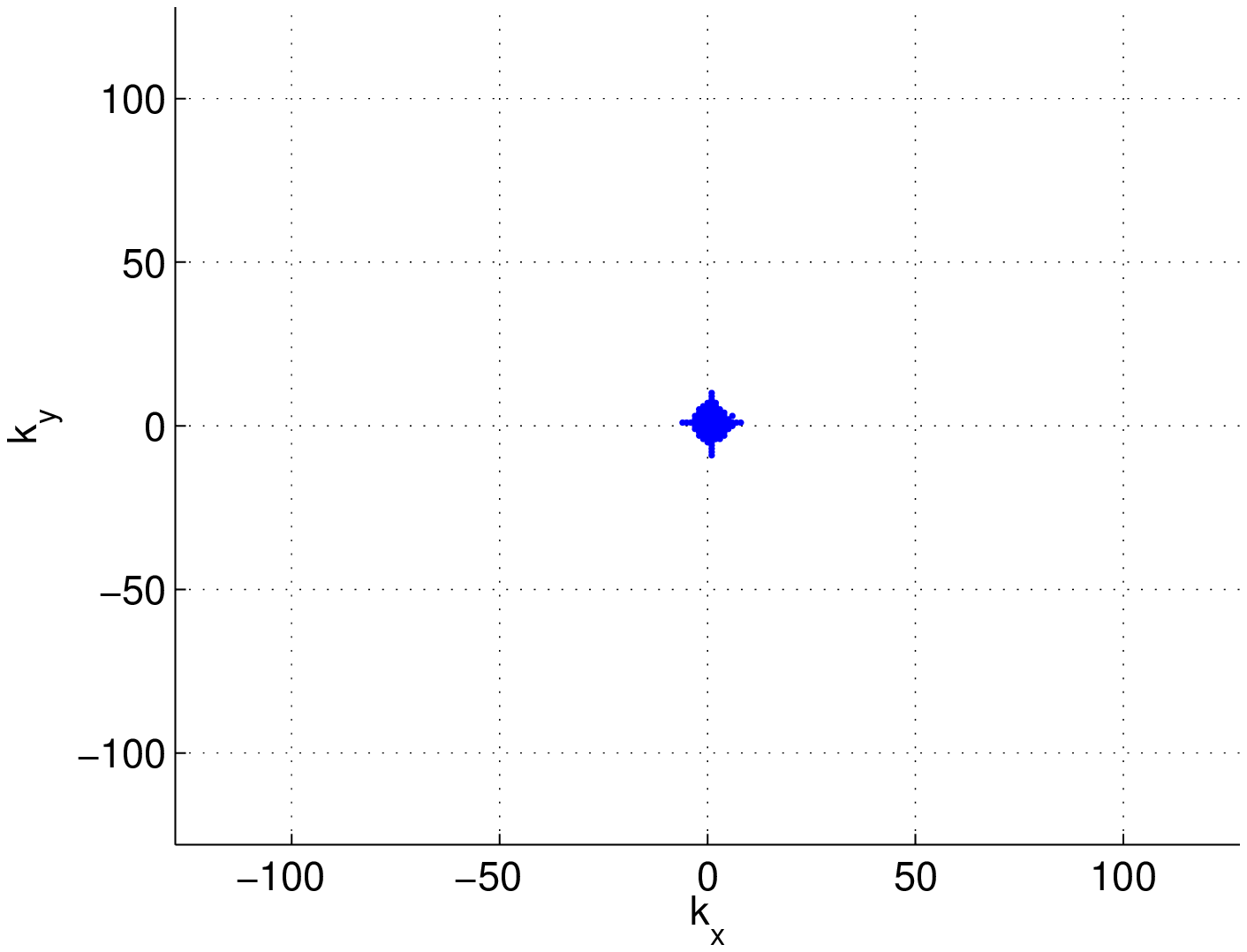}\label{phantom3dxy}}
\end{tabular}
\vspace{-0.1in}
\caption{Phantom experiment: Acquired $k-$space points used for data observation (obtained using the proposed technique) a) $3-$dimensional view, b) projected on $k_x-k_y$ plane.}\label{phantomkpoints}
\vspace{-0.1in}
\end{figure}

Obtained $k-$space points are then used in the minimum variance framework to provide posterior statistics of optical attenuation coefficient $\mu$. Similar to previous example, we considered SNR=25 while using minimum variance framework. \Fig{phantom_cvr} represents convergence of mean and variance of $\mu$ estimate versus different number of readout lines. Similar to previous example, we also performed a global optimization to find the value of $\mu$ which results in the least discrepancy between the predicted and fully sampled temperature measurement data (over the ROI). This is shown with a dashed black line in \Fig{phantom_cvr}. As one can see, mean estimate of $\mu$ converges to its optimal value by increasing the number of readout lines.
\begin{figure}[htb!!]
\centering
\vspace{-0.1in}
\includegraphics[width=3.2in]{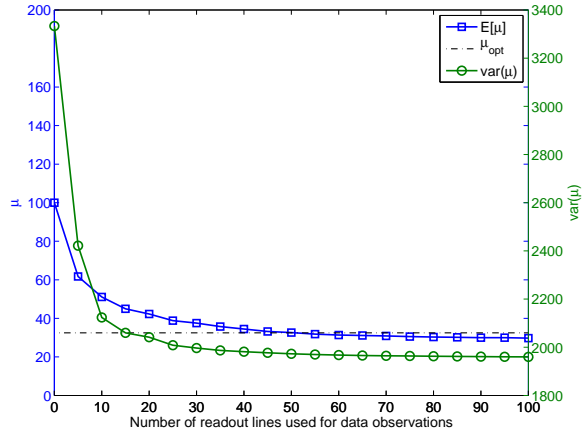}
\vspace{-0.2in}
\caption{Phantom experiment: Convergence of $\mu$ estimate versus number of observation data. Blue and dashed green lines represent mean and variance of $\mu$ estimate, respectively. Dashed black line represents the optimal value of $\mu$, obtained from global optimization.}\label{phantom_cvr}
\vspace{-0.1in}
\end{figure}
The posterior estimate of the temperature map over one of the slices, resulting via using optimally acquired $k-$space points, is shown in \Fig{tpost9}.
\begin{figure}[htb!!]
\centering
\vspace{-0in}
{\includegraphics[width=3.5in]{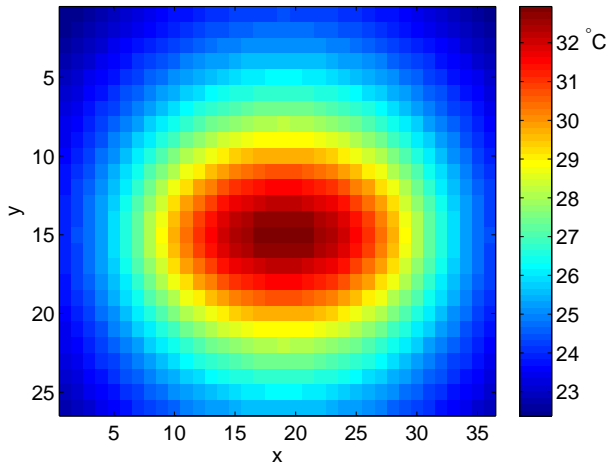}}
\vspace{-0.1in}
\caption{Phantom experiment: Posterior estimate of the temperature map over the $9^{th}$ slice (obtained using the proposed approach for $k-$space subsampling), (ROI = 18 mm $\times$ 13 mm $\times$ 26 mm). 
}\label{tpost9}
\vspace{-0.1in}
\end{figure}

\subsection*{Comparison with Other Undersampling Schemes}
We also performed model-data fusion using two other sets of data observations, which are obtained based on rectilinear and variable-density Poisson undersampling techniques. These points are shown in \Fig{phantom_unif_poiss}. Note that the same number of points, i.e. 100 readout lines, are used for data acquisition while using both schemes. In rectilinear scheme, lines are homogeneously distributed over $k_x-k_y$ plane, while most of the lines are concentrated in the center of $k_x-k_y$ plane in variable-density Poisson disk undersampling.

\begin{figure*}
\centering
\begin{tabular}{cc}
\subfigure[3-dimensional view]{\includegraphics[width=3.2in,height=2.2in]{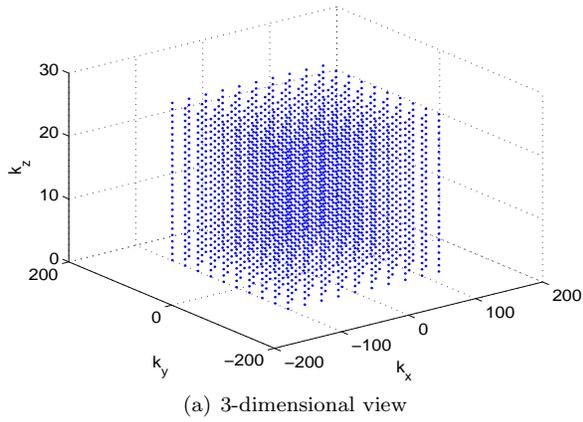}\label{phantom3dunif}} &
\subfigure[Projection over $k_x-k_y$ plane]{\includegraphics[width=3.2in,height=2.2in]{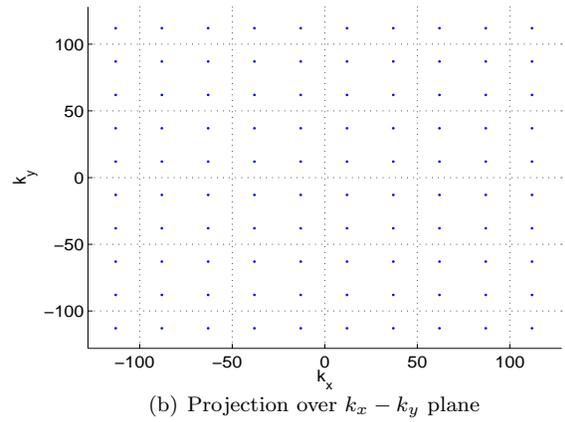}\label{phantomxyunif}} \\
\subfigure[3-dimensional view]{\includegraphics[width=3.2in,height=2.2in]{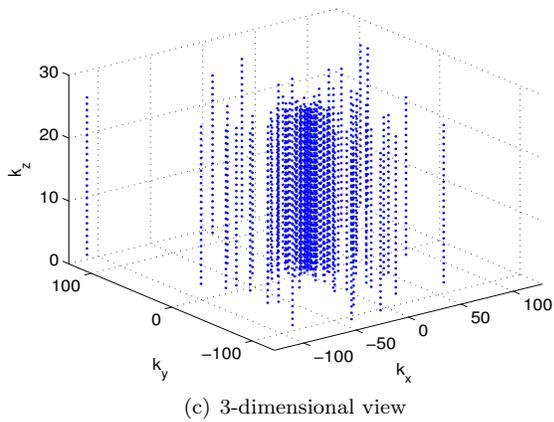}\label{phantom3dpoiss}} &
\subfigure[Projection over $k_x-k_y$ plane]{\includegraphics[width=3.2in,height=2.2in]{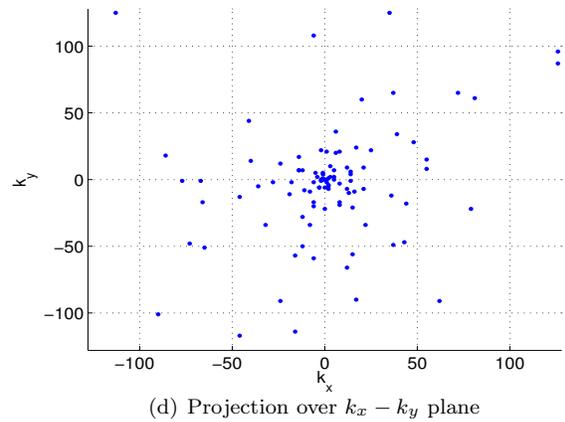}\label{phantomxypoiss}}

\end{tabular}
\vspace{-0.1in}
\caption{Phantom experiment: Acquired $k-$space points used for data observation a,b) rectilinear scheme, c,d) variable-density Poisson disk undersampling}\label{phantom_unif_poiss}
\vspace{-0.1in}
\end{figure*}

\subsection*{Error Analysis}
\Fig{phantom_err} illustrates the discrepancy between posterior temperature estimate and the temperature data acquired from fully sampled $k-$space over one of the slices, while using different data acquisition techniques. As one can see, rectilinear undersampling performs the worst (results in higher value of maximum error) among the applied techniques. Variable-density Poisson disk undersampling performs better than rectilinear scheme, but it results in slightly higher value of maximum error over the $9^{th}$ slice comparing to the proposed approach.
\begin{figure*}[htb!!]
\centering
\vspace{-0in}
\begin{tabular}{ccc}
\hspace{-0.2in}\subfigure[]{\includegraphics[width=2.5in]{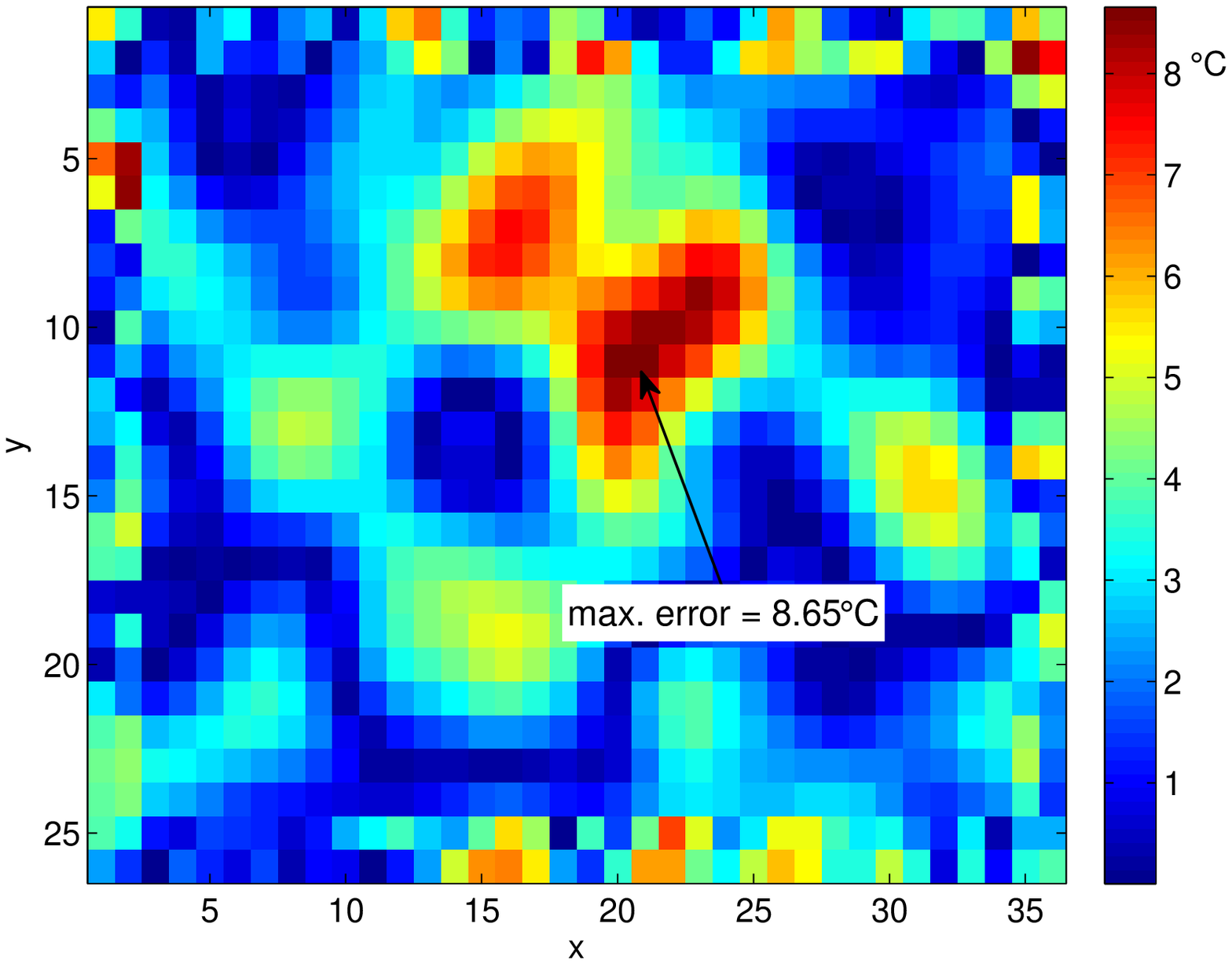}\label{errpost9}} &
\hspace{-0.2in}\subfigure[]{\includegraphics[width=2.5in]{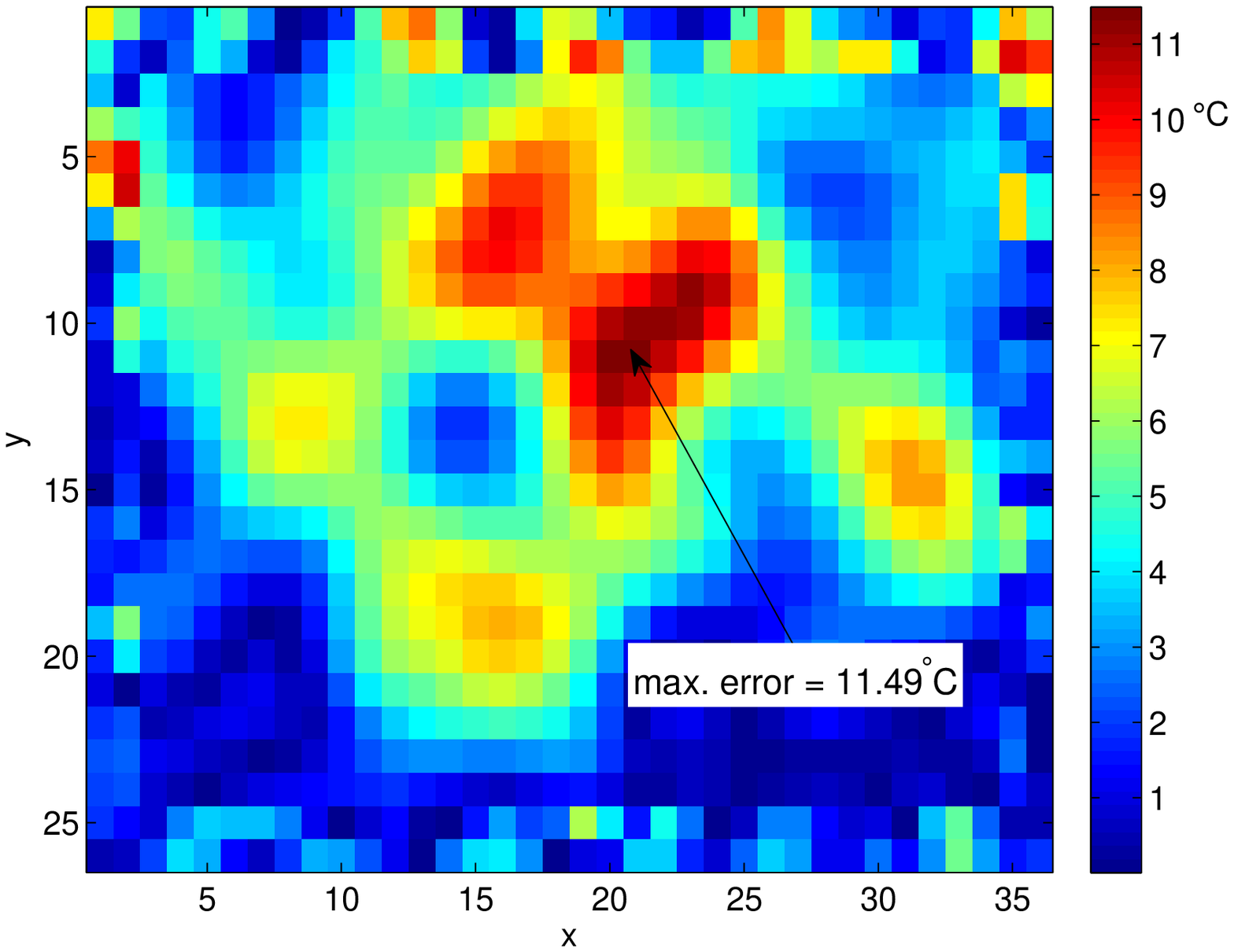}\label{err9rect}} &
\hspace{-0.2in}\subfigure[]{\includegraphics[width=2.5in]{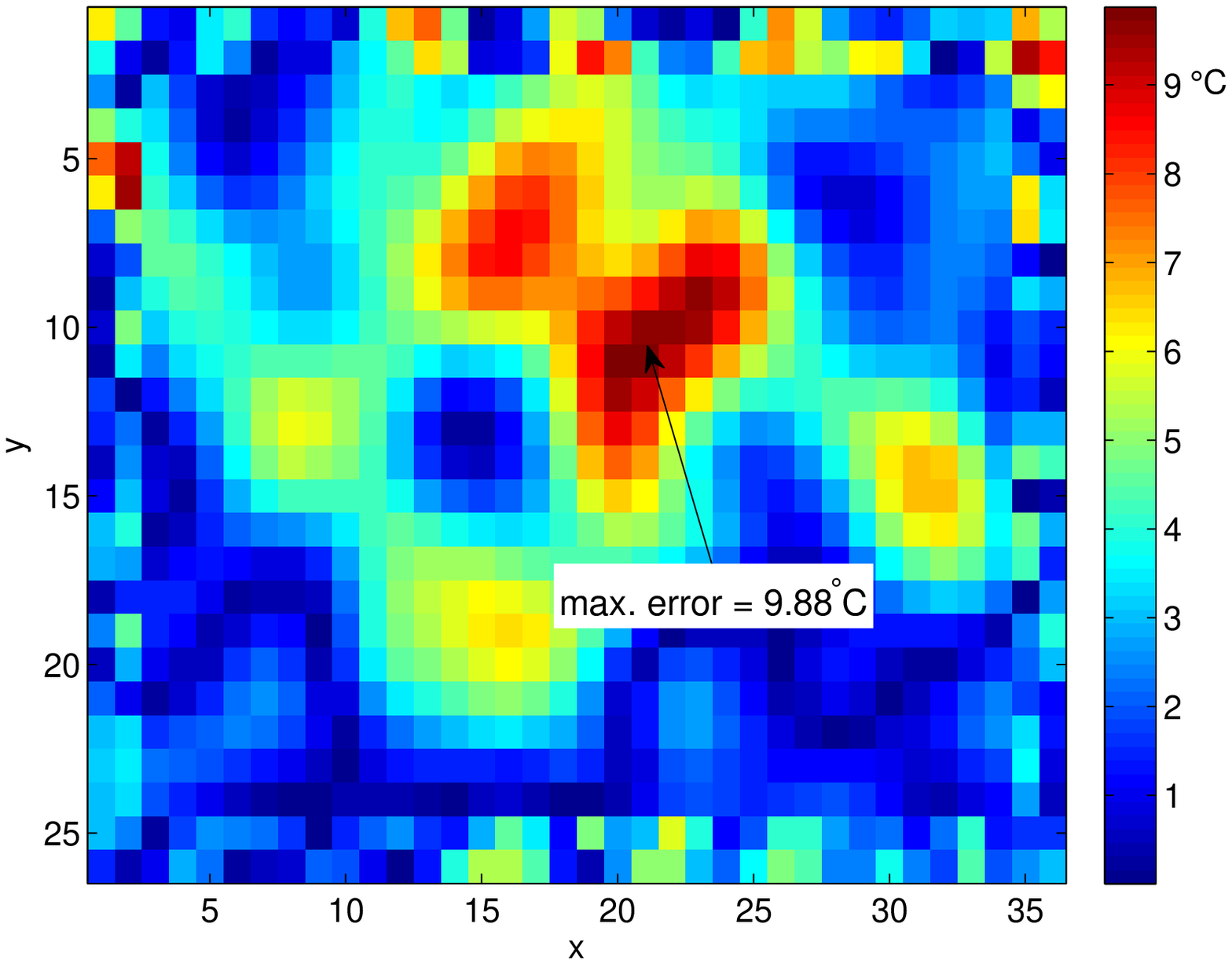}\label{err9poiss}} 
\end{tabular}
\vspace{-0.1in}
\caption{Phantom experiment: The error between the posterior estimate and actual measurement of temperature field over the $9^{th}$ slice, using a)  the proposed approach, b) rectilinear undersampling, and c) variable-density Poisson disk undersampling (ROI = 18 mm $\times$ 13 mm $\times$ 26 mm). 
}\label{phantom_err}
\vspace{-0.1in}
\end{figure*}

Similar to previous examples, we have studied the effect of increasing data observations on performance of the proposed approach. \Fig{phantom_err_cvr} represents the RMSE between the posterior temperature estimate and acquired temperature data from fully sampled $k-$space versus different numbers of readout lines used for data acquisition. As expected, RMSE reduces by increasing the number of readout lines.
\begin{figure}
\vspace{-0.1in}
\centering
\includegraphics[width=3.25in]{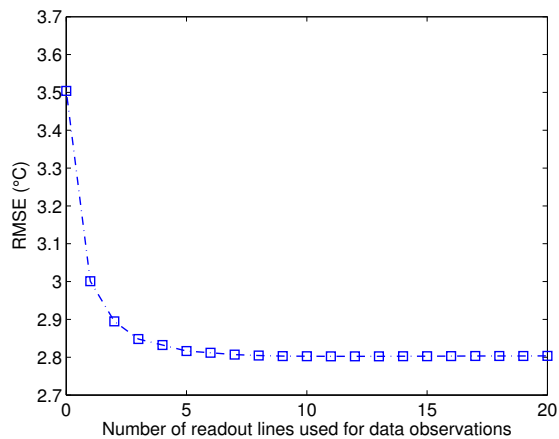}
\caption{Phantom experiment: Convergence of RMSE (over all voxels inside the ROI) between the acquired temperature map (from fully-sampled $k-$space) and temperature estimate (obtained from the proposed technique).}\label{phantom_err_cvr}
\end{figure}

Finally, Table \ref{table_phantom_err} represents the RMSE (over all voxels inside the ROI) between temperature estimates and acquired temperature data, given different number of readout lines and different data acquisition schemes. Note that there is no change in RMSE value while using rectilinear undersampling. This is expected given that none of
the acquired $k-$space points using rectilinear undersampling
covers regions with high information content (high information
region is mainly concentrated in the center of $k-$space).
Hence, almost no information is obtained by performing the
model-data fusion. Also as it can be seen in Table \ref{table_phantom_err}, proposed technique results in similar value of RMSE  comparing with variable-density Poisson disk undersampling method. 
\begin{table*}[htb!!!]
\centering
\footnotesize
\caption{Phantom experiment: Convergence of error  (over all the voxels inside the ROI) given
different numbers of $k-$space samples for model-data fusion. Prior values of error are shown in second column (with zero observations).}\label{table_phantom_err}
\begin{tabular}{lccccccc}
\hline
Number of readout lines & 0 & 5 & 10 & 30 & 50 & 75 & 100\\
\hline
RMSE (Proposed Technique) & 3.50 & 3.00 & 2.89 & 2.81 & 2.80 & 2.80 & 2.80 \\
\hline
 \multicolumn{1}{c}{RMSE (Rectilinear Undersampling)} &  3.50 & & & & & & 3.50 (using 100 readout lines) \\
\hline
 \multicolumn{1}{c}{RMSE (variable-density Poisson disk Undersampling)} &  3.50 & & & & & & 2.96 (using 100 readout lines) \\
\hline
\end{tabular}
\vspace{0.05in}
\vspace{-0.2in}
\end{table*}

{\color{black} 
\section{Discussion}\label{sec:disc}

The novelty of the proposed approach demonstrates that locations of
high-information content with respect to a model based reconstruction  of
MR thermometry may be quantitatively identified. 
The information locations identified are a consequence of the physics based
modeling of the laser induced heating and an intelligent characterization of
the measurement uncertainties of the acquisition system.
Our mathematical model of the MR thermometry acquisition is \textit{refined}
using this set of \textit{judiciously} selected data observations to
reconstruct the thermal image.
Intuitively, the accuracy and confidence of the
thermal image reconstruction depends on multiple factors like accuracy, information content,  and the 
number of data observations. 
However, in each example, the presented approach creates an accurate
assessment of temperature data while using only subset of the full
$k-$space data acquisition.  As seen in 
\Fig{mu_var_cvr},
\Fig{human2d_cvr}, and
\Fig{phantom_cvr},
the measurement locations chosen sufficiently reduce
the posterior variance of model parameters needed to create an accurate thermal image
reconstruction. 

Generally, accurate temperature image reconstructions are 
observed in the subsampled examples,
\Fig{err40},
\Fig{human2d_err}, and \Fig{phantom_err}.
Within the context of the model-data fusion, the accuracy of the Pennes bioheat
equation based temperature estimations is improved by the  measurement
locations identified by our technique. Quantiatively, 
the posterior value of RMSE is less than its prior value. 
The majority of errors in the reconstruction may be attributed to 
artifacts that are not represented in the composite sensor model.
Errors are primarily expected in
neighborhood regions of the laser fiber.


Currently, MR thermometry is acquired primarily through
quantitative phase-sensitive techniques based on the Proton Resonance
Frequency (PRF) shift. Compared with current PRF techniques,
the benefit of the proposed model-based approach is
that the predicted locations of high-information content represent a 
fraction of the fully sampled $k-$space observations.
This implies that an accurate thermal image may be reconstructed in shorter time.
Further, the typical PRF approach is susceptible to errors from
inter-/intra-scan motion, intravoxel lipid contamination, low SNR, tissue
susceptibility changes, magnetic field drift, excessive heating, and
modality$-$dependent applicator induced artifacts~\cite{poorter1995noninvasive,depoorter1994proton,kuroda1997temperature,de1999fast}.  
The potential data corruption and resulting information loss
undermines the confidence and quality in the
treatment guidance decisions provided by the MR thermometry.
The presented model-based approach provides a rigorous framework to include
a model for each source of data corruption artifact.
Improved temperature estimate accuracy in the
presence of the various sources of image contamination will facilitate a
safer, more conformal laser treatment achieved in a variety of circumstances where
critical structures may have limited its use. 

Mutual information provides the ideal measure to identify $k-$space locations
that provide the greatest possible reduction in the model parameter
uncertainty. Also, incorporating $k-$space readout trajectories will
provide the 
most physically meaningful models with respect to the underlying MR physics.
However, given non-linearities  in our physics model,
this is computationally prohibitive. In general, maximizing
the mutual information between $n$ locations over $k-$space
and our thermal image model (with $m$ unknown parameters) would result in $m+n$ dimensional integrals
to evaluate the mutual information, which can be computationally expensive due to often large values of $n$.
Here, we pursued a computationally tractable approximation and
utilized the concept of entropy maximization for approximating optimal
$k-$space readout lines.
Maximization of the mutual
information with respect to a readout trajectory in $k-$space, is a natural
progression of these efforts and is the subject of on-going work.

We emphasize that the technique presented in this manuscript is a general
methodology which can be utilized with all MR compatible thermal delivery
modalities. Algorithm development for MRgLITT is used in this manuscript as
a vehicle to focus modeling efforts and is directly extendable to heating
induced by High Intensity Focus Ultrasound (HIFU).
Highly efficient $k$-space sampling strategies for HIFU may be identified
using various models of the ultrasound based heating to 
allow reconstruction of high quality MR thermometry images with
significantly less data compared with fully sampled approaches.


\section{Conclusion}
In summary, we have presented 
a model-based information theoretic approach to improve
the efficiency of MR thermometry for monitoring MRgLITT
procedures. The approach approximates measurement locations in $k-$space
with the highest information content with respect to the thermal image
reconstruction. The performance of these predicted measurement locations 
was demonstrated
retrospectively \textit{in vivo} and \textit{in silico} imaging 
and is likely to reduce the number of $k-$space measurements needed to
accurately reconstruct the thermal image.
}

\section*{Acknowledgement}
This research work was partially supported by DoD W81XWH-14-1-0024, NIH grant TL1TR000369, and O'Donnell foundation. It was conducted at the MD Anderson Center for Advanced Biomedical Imaging in-part with equipment support from General Electric Healthcare.
\bibliographystyle{IEEEtran}
\bibliography{mda_refs,est_refs,OptimalSensor,Mypubs,Refs1,Myrefs}

%

%
%
%




\end{document}

%% file: symbols.tex
\newcommand{\re}{\mathbb R}
\newcommand{\Ex}[2]{\mathcal{E}^{#2}[#1]}

\newcommand{\etal}{\textit{et al. }}

\newcommand{\eq}[1]{Eq.~(\ref{#1})}
\newcommand{\Fig}[1]{Fig.~\ref{#1}}

\newcommand{\xii}{\boldsymbol{\xi}}

\newcommand{\x}{\mathbf{x}}

\newcommand{\z}{\mathbf{z}}